\documentclass[times,sort&compress,3p]{elsarticle}
\journal{Journal of Multivariate Analysis}
\usepackage[labelfont=bf]{caption}

\usepackage{caption,subcaption,mwe}

\usepackage{amsmath,amsfonts,amssymb,amsthm,booktabs,color,epsfig,graphicx,hyperref,url}

\theoremstyle{plain}% Theorem-like structures provided by amsthm.sty

\newtheorem{rem}{Remark}
\newtheorem{proposition}{Proposition}
\newtheorem{lemma}{Lemma}

\theoremstyle{definition}

\newcommand{\bb}[1]{\boldsymbol{#1}}
\newcommand{\e}{\varepsilon}
\newcommand{\EE}{\mathsf{E}}
\newcommand{\ind}{\mathbf{1}}
\newcommand{\N}{\mathbb{N}}
\newcommand{\oo}{\mathrm{o}}
\newcommand{\OO}{\mathcal O}
\newcommand{\PP}{\mathsf{P}}
\newcommand{\R}{\mathbb{R}}
\newcommand{\rd}{d}
\newcommand{\VV}{\mathsf{var}}

\begin{document}

\begin{frontmatter}

\title{Minimax properties of Dirichlet kernel density estimators}

    \author[a1]{Karine Bertin\texorpdfstring{}{)}}
    \author[a2]{Christian Genest\corref{mycorrespondingauthor}}
    \author[a3]{Nicolas Klutchnikoff}
    %\ead{nicolas.klutchnikoff@univ-rennes2.fr}
    \author[a2,a4]{Fr\'ed\'eric Ouimet}
    %\ead{frederic.ouimet@umontreal.ca}

    \address[a1]{CIMFAV-INGEMAT, Universidad de Valpara\'iso, General Cruz 222, Valpara\'iso, Chile}%
    \address[a2]{Department of Mathematics and Statistics, McGill University, Montr\'eal (Qu\'ebec) Canada H3A 0B9}
    \address[a3]{Univ Rennes, CNRS, IRMAR - UMR 6625, F-35000 Rennes, France}%
    \address[a4]{Centre de recherches math\'ematiques, Universit\'e de Montr\'eal, Montr\'eal (Qu\'ebec) Canada H3T 1J4}

\cortext[mycorrespondingauthor]{Corresponding author. Email address: christian.genest@mcgill.ca \url{https://www.math.mcgill.ca/cgenest/}}

\begin{abstract}
This paper considers the asymptotic behavior in $\beta$-H\"older spaces, and under $L^p$ losses, of a Dirichlet kernel density estimator proposed by \citet{doi:10.2307/2347365} for the analysis of compositional data. In recent work, \citet{MR4319409} established the uniform strong consistency and asymptotic normality of this estimator. As a complement, it is shown here that the Aitchison--Lauder estimator can achieve the minimax rate asymptotically for a suitable choice of bandwidth whenever $(p,\beta) \in [1, 3) \times (0, 2]$ or $(p, \beta) \in \mathcal{A}_d$, where $\mathcal{A}_d$ is a specific subset of $[3, 4) \times (0, 2]$ that depends on the dimension $d$ of the Dirichlet kernel. It is also shown that this estimator cannot be minimax when either $p \in [4, \infty)$ or $\beta \in (2, \infty)$. These results extend to the multivariate case, and also rectify in a minor way, earlier findings of \citet{MR2775207} concerning the minimax properties of Beta kernel estimators.
\end{abstract}


\end{frontmatter}

\section{Introduction\label{sec:1}}

Compositional data refer to observations of a random vector whose components represent proportions of a whole whose size is either irrelevant or analyzed separately; see, e.g., the book by Aitchison~\cite{MR865647} for an introduction to this topic. As proportions are non-negative and sum up to $1$ by definition, one can write a compositional vector of length $d + 1$ in the form $(\bb{X}, X_{d+1}) = (\bb{X}, 1 - \|\bb{X}\|_1)$ with $\bb{X}$ belonging to the $d$-dimensional simplex
\[
\mathcal{S}_d = \big\{\bb{s} \in [0,1]^d : \| \bb{s} \|_1 \leq 1 \big\},
\]
where $\|\bb{s}\|_1 = |s_1| + \cdots + |s_d|$ denotes the $\ell^1$ norm on $\mathbb{R}^d$. When all proportions are known to be strictly positive, the vector $\bb{X}$ belongs to the interior of $\mathcal{S}_d$, denoted $\mbox{Int} (\mathcal{S}_d)$.

As illustrated, e.g., by Filzmoser et al.~in their book~\cite{MR3839314}, compositional data arise in a wide range of fields such as chemometrics, demography, economics, geochemistry, and survey methodology. In dimension $d + 1 \ge 3$, the first and best known approach to modeling compositional data is to transform them to the $d$-dimensional simplex through an additive log-ratio map. This strategy can be used to construct a kernel on $\mbox{Int} (\mathcal{S}_d)$ by applying the logistic transformation to the multivariate Gaussian density; see, e.g., \cite{doi:10.2307/2347365, doi:10.1016/j.cageo.2009.12.011}. Possible alternatives include boundary kernels derived as a solution to a variational problem \citep{MR1680306} and products of one-dimensional asymmetric kernels \citep{MR2568128}; see also \cite{MR3760293, Kokonendji_Some_2021} for a general theory of multivariate associated kernels. For a survey of asymmetric kernel methods, including their use for compositional data in arbitrary dimension, refer to \cite{MR4319409}.

Difficulties with the transformation approach may occur when the data are suspected to be sparse or include zeros, either for structural reasons, rounding or otherwise \cite{doi:10.1023/A:1023866030544}. An alternative initially due to~\citet{doi:10.2307/2347365} and later considered by \citet{doi:10.1016/j.cageo.2009.12.011} is to work directly on the simplex. Their strategy consists of basing a kernel density estimator on the Dirichlet distribution with arbitrary parameters $\bb{u} = (u_1, \dots, u_d) \in (0, \infty)^d$ and $v \in (0,\infty)$, whose density is given, for all $\bb{s} \in \mathcal{S}_d$,~by
\begin{equation}
\label{eq:1}
K_{\bb{u},v}(\bb{s}) = \frac{\Gamma(v + u_1 + \cdots + u_d)}{\Gamma (v) \Gamma (u_1) \cdots \Gamma (u_d)} \, (1 - \|\bb{s}\|_1)^{v - 1} \prod_{i=1}^d s_i^{u_i - 1}.
\end{equation}

Given a random sample $\bb{X}_1, \ldots,\bb{X}_n$ from an unknown density $f$ on $\mathcal{S}_d$, \citet{doi:10.2307/2347365} suggest that a Dirichlet kernel density estimator with bandwidth parameter $b \in (0, \infty)$ could then be defined, for all $\bb{s}\in \mathcal{S}_d$, by
\begin{equation}
\label{eq:2}
f_{n,b}(\bb{s}) = \frac{1}{n} \sum_{i=1}^n K_{{\bb{s}}/{b} + \bb{1}, {(1 - \|\bb{s}\|_1)}/{b} + 1}(\bb{X}_i),
\end{equation}
where $\bb{1} = (1, \ldots, 1)$ is a $d$-dimensional vector whose components are all equal to $1$.

When $d = 1$, the Dirichlet density in~\eqref{eq:1} coincides with the Beta density defined, for all $s \in (0, 1)$, by $K_{u,v} (s) = \Gamma(u + v) (1 - s)^v s^u / \{\Gamma(u) \Gamma(v)\}$. The associated estimator in~\eqref{eq:2} then corresponds to the Beta kernel density estimator whose theoretical properties were initially investigated by \citet{MR1685301} and \citet{MR1718494}. For additional work on this topic, see, e.g., \cite{MR2775207, MR3333996}, \cite{MR1985506, MR2568128, MR4148613}, \cite{Charpentier2006phd, Charpentier_Fermanian_Scaillet_2007}, \cite{MR2206532}, \cite{MR2756441, MR4471289, MR3463548, MR4076244, MR4136585}, \cite{doi:10.1016/j.jbankfin.2003.10.018}, \cite{MR2598955}.

A key feature of the Aitchison--Lauder proposal is that the shape of the kernel changes with the position $\bb{s}$ in the simplex. This makes it possible to avoid the boundary bias problem associated with traditional estimators in which the kernel is the same at every point. An alternative strategy in which the bandwidth parameter $b$ is exclusively related to the concentration of the distribution was proposed by \citet{Martin/etal:2006}.

To this point, the asymptotic properties of kernel-based estimators defined directly on the simplex have been limited to the studies of \citet{MR1293514}, \citet{MR4287788}, and \citet{MR4319409}. The first two papers were concerned with the properties of Bernstein estimators in dimension $d + 1 = 3$ and above, respectively. In contrast, \citet{MR4319409} focused on the Aitchison--Lauder estimator, whose uniform strong consistency and asymptotic normality were established. In particular, these authors showed that in any dimension, the estimator defined in~\eqref{eq:2} achieves the optimal convergence rate $\mathcal{O} (n^{-4/(d+4)})$ for the mean squared error and the mean integrated squared error when the underlying density is twice continuously differentiable on~$\mathcal{S}_d$.

In this paper, the performance of the Dirichlet kernel estimator defined in~\eqref{eq:2} is studied from the point of view of asymptotic minimax theory. To this end, it will be assumed that the unknown density $f$ is sufficiently smooth that, for some regularity parameter $\beta  \in (0, \infty)$ and Lipschitz constant $L \in (0, \infty)$, $f$ belongs to the H\"older space
\[
\Sigma(d, \beta, L) = \left\{ f: \mathcal{S}_d \to \mathbb{R} ~:~\forall_{\bb{\gamma} \in \mathbb{N}_0^d: \|\bb{\gamma}\|_1 = m} \;
\forall_{\bb{s},\bb{t} \in \mathrm{Int}(\mathcal{S}_d)} \; | D^{\bb{\gamma}} f(\bb{s}) - D^{\bb{\gamma} }f(\bb{t})|  \leq L \, \|\bb{s} - \bb{t}\|_1^{\beta - m} \right\},
\]
where $\mathbb{N}_0 = \{ 0, 1, \ldots \}$, $m = \sup \{\ell\in \mathbb{N}_0 : \ell < \beta \}$, and for every $\bb{\gamma} = (\gamma_1, \ldots, \gamma_d) \in \mathbb{N}_0^d$ and $\bb{s} = (s_1, \ldots, s_d) \in \mathcal{S}_d$,
\[
D^{\bb{\gamma}} f (\bb{s}) = \frac{\partial^{\bb{\gamma}}}{\partial s_1^{\gamma_1} \cdots \partial s_d^{\gamma_d}} \, f (\bb{s})
\]
with the convention that if $\bb{\gamma} = \bb{0}$, the $d$-dimensional vector whose components are all equal to $0$, then $D^{\bb{\gamma}} f = f$.

Of interest here is whether or not the Aitchison--Lauder estimator can possibly achieve the minimax rate of convergence on $\Sigma(d, \beta, L)$ for any given $L^p$ loss defined, for any $p \in [1, \infty)$ and estimator $f_n$ of $f$, by
\[
R_n (f_n, f) = \big \{ \EE \big( \| f_n - f\|_p^p \big) \big\}^{1/p}
\]
whenever this expectation exists. The corresponding risk of the estimator $f_n$ over the class $\Sigma(d, \beta, L)$ is then given by
\[
R_n (f_n, d, \beta, L) = \sup_{f\in \Sigma(d, \beta, L)} R_n (f_n, f)
\]
and the minimax rate of convergence over $\Sigma(d, \beta, L)$ is defined by
\[
r_n (d, \beta, L) = \inf_{f_n} R_n (f_n, d, \beta, L),
\]
where the infimum is taken over all possible estimators $f_n$ of $f$. From Theorem~2 and Remark~3 of~\citet{MR4029144}, this rate is known to be
\begin{equation}
\label{eq:3}
r_n (d, \beta, L) \asymp_{d, \beta, L} \varphi_n (d, \beta) = n^{-\beta /(d + 2 \beta)}.
\end{equation}

It is first shown in Section~\ref{sec:2} that for any $L^p$ loss with $p \in [1, 3)$ and density $f \in \Sigma(d, \beta, L)$ with regularity parameter $\beta \in (0, 2]$, there exists a sequence of Dirichlet kernel density estimators with suitably chosen bandwidth parameter which achieves the minimax rate~\eqref{eq:3} asymptotically. This result is also proved for pairs $(p,\beta)$ in a specific subset of $[3,4) \times (0,2]$ denoted by $\mathcal{A}_d$ and defined in~\eqref{eq:4} in Section~\ref{sec:2}. As detailed in Proposition~\ref{prop:1}, this bandwidth parameter depends both on the sample size, $n \in \mathbb{N} = \{ 1, 2, \ldots \}$, and on the regularity parameter $\beta$ of the underlying density~$f$. As the value of $\beta$ is typically unknown, this result is primarily of theoretical interest but could motivate the search for data-driven bandwidth selection procedures based on cross-validation or Goldenshluger--Lepski-type procedures in this context; see Remark~\ref{remark:2} below.

As shown in Section~\ref{sec:3}, however, the Aitchison--Lauder class of kernel estimators cannot achieve the asymptotic minimax rate for densities having a high degree of smoothness, namely $\beta \in (2, \infty)$, or if the reference loss function is $L^p$ for some $p \in [4, \infty)$ when $\beta\in (0,2]$. These results constitute Propositions~\ref{prop:2} and~\ref{prop:3}, respectively. Alas, the techniques used to prove these results are inadequate to settle the case $(p, \beta) \in [3, 4) \times (0, 2] \setminus \mathcal{A}_d$, which remains open.

The results presented here extend to all dimensions previous findings of \citet{MR2775207} in dimension $d + 1 = 2$. These authors showed that for a suitable choice of bandwidth, Beta kernel density estimators of an unknown density $f \in \Sigma(d, \beta, L)$ can achieve the minimax rate asymptotically for any $L^p$ loss with $p \in [1, 3)$ when $\beta \in (0, 2]$, but not when $(p, \beta) \in [2, \infty) \times (2, \infty)$ or $(p, \beta) \in [4, \infty) \times (0, 2]$. The cases $(p, \beta) \in [3, 4) \times (0, 2]$ and $(p, \beta) \in [1, 2) \times (2, \infty)$ were not covered in \cite{MR2775207}. An oversight in one of the arguments presented in that paper is corrected along the way.

In Section~\ref{sec:4}, insight into the meaning of the regularity parameter $\beta$ is provided by determining its value for the classical Dirichlet distributions, and concluding comments are given in Section~\ref{sec:5}. Throughout the paper, expectation is taken with respect to the joint law of the mutually independent copies $\bb{X}_1, \ldots, \bb{X}_n$ of $\bb{X}$. Whether explicitly or not, the bandwidth parameter $b$ is always assumed to be a function of the sample size except in Lemmas~\ref{lemma:B2}~and~\ref{lemma:C3} and their proofs. The notation $u = \OO(v)$ means that $\limsup |u / v| < C < \infty$ as $n \to \infty$ or $b \to 0$, where the strictly positive constant $C$ depends on no variable unless explicitly written as a subscript. In some instances, $u \ll v$ is used to mean $u \geq 0$ and $u = \OO(v)$. When $u \ll v$ and $v \ll u$ both hold, then one writes $u \asymp v$. Similarly, the notation $u = \oo(v)$ means that $\lim |u / v| = 0$ as $n\to \infty$ or $b\to 0$. Subscripts indicate which parameters the convergence rate can depend on.

\section{Minimax result\label{sec:2}}

The following result identifies a set of $L^p$ loss functions and a range of regularity parameter values $\beta$ on the unknown density $f \in \Sigma(d, \beta, L)$ for which the minimax rate can be achieved by the Aitchison--Lauder estimator defined in~\eqref{eq:2} with a suitable choice of bandwidth parameter $b \in (0, \infty)$. This choice of bandwidth depends both on the sample size $n$ and on the degree $\beta$ of smoothness of the underlying density $f$.

\begin{proposition}
\label{prop:1}
Assume that $(p,\beta) \in [1,3) \times (0,2]$ or that $(p,\beta) \in \mathcal{A}_d$, where
\begin{equation}
\label{eq:4}
\mathcal{A}_d = \left\{(p,\beta)\in [3,4) \times (0,2] : 3 \leq p \leq 2 + \frac{d}{d - 1}, \, 2d \, \frac{p-3}{p-2} < \beta \leq 2\right\}.
\end{equation}
Moreover, let $b_n = c n^{-2/(d + 2 \beta)}$ for every integer $n \in \mathbb{N}$ and some arbitrary constant $c \in (0, \infty)$. Then the sequence $\{f_{n,b_n} : n \in \mathbb{N} \}$ achieves the minimax rate, namely
\[
\limsup_{n\to \infty} \frac {R_n(f_{n,b_n}, d, \beta, L)} {r_n(d, \beta, L)} < \infty.
\]
\end{proposition}

\bigskip
The proof is deferred to~\ref{app:A}. The following comments are in order.

\begin{rem}
\label{remark:1}
\emph{In dimension $d +1 = 2$, the statement of Proposition~\ref{prop:1} extends Theorem~1 of~\citet{MR2775207} from $[1,3) \times (0,2]$ to the set $\mathcal{A}_d$. In that paper, however, the range of $p$ was mistakenly claimed to be the interval $[1, 4)$ due to a slip near the end of the proof. Indeed, upon taking into account the integrability conditions
\[
0 \leq p/4 +  \max (0, p - 2)/4 < 1, \quad 0 \leq p/4 < 1
\]
that appear in \cite{MR2775207}, Theorem~1 therein is only true for $(p, \beta) \in [1, 3) \times (0, 2]$ rather than for $(p, \beta) \in [1, 4) \times (0, 2]$.}

\emph{Given that the union $[1,3) \times (0,2] \cup \mathcal{A}_d$ is strictly included in $[1, 4) \times (0, 2]$ for every integer $d \in \N$, the reader might wonder what is the difficulty in extending the result of Proposition~\ref{prop:1} to all $(p, \beta) \in [1, 4) \times (0, 2]$. When $p \in (0, 2]$, the centered absolute $p$th moment of the estimator $f_{n,b}(\bb{s})$ can be bounded above using Jensen's inequality, viz.
\[
\EE \big[ |f_{n,b}(\bb{s}) - \EE \{ f_{n,b}(\bb{s}) \} |^p \big] \leq \big[ \VV \big\{ f_{n,b}(\bb{s}) \big\} \big]^{p/2}.
\]
The variance term is then relatively easy to control using the asymptotics of the gamma function.}

\emph{If this upper bound were to hold for all $p \in (0,4)$, one could then extend the result of Proposition~\ref{prop:1} to all $(p, \beta) \in [1, 4) \times (0, 2]$ by using the same control on the variance term. However, when $p \in (2,\infty)$, this bound is no longer valid because the function $x\mapsto x^{2/p}$ is now concave instead of convex.
It has to be adjusted to take into account the supremum norm $\|\cdot\|_{\infty}$ of the summands in $f_{n,b}(\bb{s})$ as follows:
\[
\EE \big[ |f_{n,b}(\bb{s}) - \EE \{ f_{n,b}(\bb{s}) \} |^p \big] \ll_p \left\{\|K_{{\bb{s}}/{b} + \bb{1}, {(1 - \|\bb{s}\|_1)}/{b} + 1} \|_{\infty}/n\right\}^{p-2} \VV \big\{ f_{n,b}(\bb{s}) \big\} + \big[ \VV \big\{ f_{n,b}(\bb{s}) \big\} \big]^{p/2}.
\]
See~(4.11) of \citet{MR523165}.}

\emph{Upon applying the best local bounds available on both the supremum and variance terms, and then integrating on both sides, see~\eqref{eq:A.6}, one can see that the supremum term will impose the slightly more stringent integrability condition $p < 3$, but the minimax rate will still be reached for all $(p,\beta) \in [1,3) \times (0,2]$. By applying a log-convex mixture between the local and uniform bounds on $\|K_{\bb{s}/b + \bb{1}, (1 - \|\bb{s}\|_1)/b + 1}\|_{\infty}$ in the above reasoning, see~\eqref{eq:A.5}, one can then extend the minimax results to also include all $(p,\beta)\in \mathcal{A}_d$. It remains a mystery whether another argument could push the boundary of the minimax results even further inside $[1, 4) \times (0, 2]$.}
\end{rem}

\begin{rem}
\label{remark:2}
\emph{When attention is restricted to the $L^2$ loss, Proposition~\ref{prop:1} implies that Dirichlet estimators achieve the minimax rate of convergence when $\beta$ is known and belongs to $(0,2]$. Adaptive estimators, which converge at the minimax rate irrespective of the smoothness $\beta$, may be obtained by defining a statistical procedure akin to the one developed by Goldenshluger and Lepski; see \cite{MR2850214,MR3230001} and references therein. The main idea behind this method is to select a data-driven bandwidth $\hat{b}$, in a finite set $\mathcal{B}$ ranging from  $n^{-1/d}$ to $1$, that satisfies
\begin{equation}
\label{eq:5}
\hat{b}=\min_{b\in\mathcal{B}}   \left\{  A(b)+V(b)    \right\},
\end{equation}
where, for any $\delta \in (0, \infty)$,
\[
V(b) = {F}_\infty C_1(1 + \delta) / (nb^{d/2}) > \EE \big\{ \|f_{n, b} - \EE ( f_{n, b}) \|_2^2 \big\}.
\]
The latter quantity is a penalized version of the integrated variance of the estimator $f_{n, b}$ while
\[
A(b) = \max_{b' \in \mathcal{B}} \big\{ \|f_{n, b} - f_{n, b \vee b' }\|_2^2 - V(b) - V (b \vee b') \big\}_+
\]
is an estimation of its integrated bias term. Here, the notations $a \vee b = \max (a, b)$ and $a_+ = a \vee 0$ are used, the density $f$ is assumed to satisfy $\|f\|_{\infty} \le {F}_{\infty}$, and $C_1 = 2^{-d} \sqrt{\pi} / \Gamma(d/2 + 1/2)$ following~(15) and~(17) of \citet{MR4319409}, together with the integral calculation in~(4.10) of \citet{MR3825458}.}

\emph{Thus,~\eqref{eq:5} can be interpreted as an empirical version of the usual bias-variance trade-off. The study of this procedure relies on the fine control of the process $\|f_{n,b} - \EE(f_{n,b})\|_2$ through the probabilities $\PP \{ \|f_{n,b} - \EE ( f_{n,b}) \|_2 > t \}$, which can be bounded above using Bernstein's inequality for $U$-statistics and Hoeffding's inequality; see \cite{MR3333996} for more details. Using a similar approach, one can prove that the resulting estimator satisfies
\[
\limsup_{n\to\infty} \frac{R_n(f_{n,\hat{b}}, d, \beta, L)}{r_n(d,\beta,L)} < \infty
\]
for any $\beta\in(\min\left(2, d/2\right), 2]$. Observe that this adaptive result only makes sense for $d \in \{ 1, 2, 3 \}$ and that the range of values of $\beta$ for which it holds is rather small, especially when $d = 3$.}
\end{rem}

\section{Non-minimaxity results\label{sec:3}}

The following results identify a set of $L^p$ loss functions and a range of regularity parameter values $\beta$ on the unknown density $f \in \Sigma(d, \beta, L)$ for which the Aitchison--Lauder estimator defined in~\eqref{eq:2} cannot possibly achieve the minimax rate of convergence, irrespective of the choice of bandwidth parameter $b \in (0, \infty)$. The values $(p, \beta) \in [1, 2) \times (2, \infty)$ were not covered in the case $d = 1$ studied by~\citet{MR2775207}.

\begin{proposition}
\label{prop:2}
Let $p \in [1, \infty)$ and $\beta \in (2, \infty)$. Then, for all sequences $\{b_n : n \in \mathbb{N} \}$ in $(0,1)$, the family $\{ f_{n,b_n}: n \in \mathbb{N} \}$ of estimators satisfies
\[
\liminf_{n\to \infty} \frac {R_n(f_{n,b_n}, d, \beta, L)} {r_n(d, \beta, L)} = \infty.
\]
\end{proposition}

\bigskip
\begin{proposition}
\label{prop:3}
Let $p \in [4, \infty)$ and $\beta \in (0, 2]$. Then, for all sequences $\{b_n : n \in \mathbb{N} \}$ in $(0,1)$, the family $\{f_{n,b_n} : n \in \mathbb{N} \}$ of estimators satisfies
\[
\liminf_{n\to \infty} \frac {R_n(f_{n,b_n}, d, \beta, L)} {r_n(d, \beta, L)} = \infty.
\]
\end{proposition}
The proofs are deferred to~\ref{app:B}.

\section{Regularity of the Dirichlet distributions\label{sec:4}}

The critical technical condition under which Propositions~\ref{prop:1}--\ref{prop:3} are established is the assumption that the underlying $d$-variate density $f$ belongs to the $\beta$-H\"older space $\Sigma (d, \beta, L)$ for appropriate choices of regularity parameter $\beta \in (0, \infty)$ and Lipschitz constant $L \in(0,\infty)$. To get a better feel for this requirement --- and the conditions on $\beta$ which guarantee that the Aitchison--Lauder estimator can achieve the asymptotic minimax rate for a suitable choice of bandwidth parameter  --- the case where $f$ is a Dirichlet density is briefly considered in this section.

Let $f = K_{\bb{u},v}$ be a $d$-variate Dirichlet density of the form~\eqref{eq:1} for some choice of parameters $\bb{u} = (u_1, \ldots, u_d) \in (0, \infty)^d$ and $v \in (0,\infty)$. First note that if $u_1, \ldots, u_d$ and $v$ are all positive integers, then $f$ is a finite multivariate polynomial. Consequently, it belongs to the H\"older space $\Sigma (d, \beta, L)$ for every $\beta \in (0, \infty)$, i.e.,
\[
\forall_{\beta \in (0, \infty)} \; \exists_{L \in (0, \infty)} \quad K_{\bb{u},v} \in  \Sigma (d, \beta, L).
\]

Next, assume that $u_1, \ldots, u_d, v \in (1, \infty)$ but that not all of them are integers. Let $m = \sup \{\ell\in \mathbb{N}_0 : \ell < \beta \}$ with
\[
\beta = \min(u_1 - 1, \dots, u_d - 1, v - 1).
\]
Then $K_{\bb{u},v} \in \Sigma(d, \beta, L)$ for some large enough Lipschitz constant $L = L (d, u_1, \ldots, u_d, v) \in (0, \infty)$. Indeed, for each integer $j \in \{ 1, \ldots, d\}$, one can write the $j$th partial derivative of $K_{\bb{u},v}(\bb{s})$ with respect to $\bb{s}\in \mathrm{Int}(\mathcal{S}_d)$ as
\[
\frac{\partial}{\partial s_j} \, K_{\bb{u},v}(\bb{s}) = (v + u_1 + \dots + u_d - 1) \, \{K_{\bb{u} - \bb{e}_j, v}(\bb{s}) - K_{\bb{u}, v - 1}(\bb{s}) \},
\]
where $\bb{e}_j$ denotes the $j$th standard basis vector in $\R^d$. Hence, for any $\bb{\gamma}\in \mathbb{N}_0^d$ such that $\|\bb{\gamma}\|_1 = m\in \N_0$, one has
\[
D^{\bb{\gamma}} K_{\bb{u},v}(\bb{s}) = \left\{\prod_{i=1}^m (v + u_1 + \dots + u_d - i)\right\} \sum_{j_1=0}^{\gamma_1} \dots \sum_{j_d=0}^{\gamma_d} \binom{\gamma_1}{j_1} \dots \binom{\gamma_d}{j_d} (-1)^{m - \|\bb{j}\|_1} K_{\bb{u} - \bb{j}, v - (m - \|\bb{j}\|_1)}(\bb{s}),
\]
with the convention that a product over an empty set equals $1$ in the trivial case $m = 0$.

Observe that the factors $K_{\bb{u} - \bb{j}, v - (m - \|\bb{j}\|_1)}(\bb{s})$ only consist of constants depending on $u_1,\dots,u_d$ and $v$ multiplied by monomials in the variables $s_1, \ldots, s_d$ and $s_{d+1} = 1 - \|\bb{s}\|_1$. Therefore, one can conclude to the existence of a constant $C \in (0, \infty)$ depending only on the integer $d$ and the reals $u_1, \ldots, u_d$, and $v = u_{d+1}$ such that, letting $j_{d+1} = m - \|\bb{j}\|_1$,
\[
|D^{\bb{\gamma}} K_{\bb{u},v}(\bb{s}) - D^{\bb{\gamma}} K_{\bb{u},v}(\bb{t})| \leq C \sum_{j_1=0}^{\gamma_1} \dots \sum_{j_d=0}^{\gamma_d} \left|\prod_{i=1}^{d+1} s_i^{u_i - 1 - j_i} - \prod_{i=1}^{d+1} t_i^{u_i - 1 - j_i}\right| .
\]

Now for arbitrary reals $p_1, \dots, p_{d+1} \in (0, \infty)$, set $p_* = \min(p_1, \ldots, p_{d+1},1)$ and $p^* = \max(p_1, \ldots, p_{d+1},1)$. One can apply a simple chaining argument to show that, for any integer $d \in \N$, all $(x_{1,1}, \ldots, x_{1,d+1}), (x_{2,1}, \ldots, x_{2,d+1}) \in [0,1]^{d+1}$, and every real $p_1, \dots, p_{d+1} \in (0, \infty)$, one has
\begin{align*}
\left|\prod_{i=1}^{d+1} x_{1,i}^{p_i} - \prod_{i=1}^{d+1} x_{2,i}^{p_i}\right|
&= \left|\sum_{k=1}^{d + 1} \left\{\left(\prod_{i=1}^{k-1} x_{1,i}^{p_i}\right) \left(\prod_{i=k}^{d+1} x_{2,i}^{p_i}\right) - \left(\prod_{i=1}^k x_{1,i}^{p_i}\right) \left(\prod_{i=k+1}^{d+1} x_{2,i}^{p_i}\right)\right\}\right| \\[2mm]
&\leq \sum_{k=1}^{d+1} |x_{1,k}^{p_k} - x_{2,k}^{p_k}| \leq \sum_{k=1}^{d+1} \left\{|x_{1,k} - x_{2,k}|^{p_k} \ind_{(0,1)}(p_k) + |p_k| \times |x_{1,k} - x_{2,k}| \ind_{[1,\infty)}(p_k)\right\} \\
&\leq p^* \sum_{k=1}^{d+1} |x_{1,k} - x_{2,k}|^{p_*} \leq p^*  (d + 1) \left(\sum_{k=1}^{d+1} |x_{1,k} - x_{2,k}|\right)^{p_*}.
\end{align*}
Applying this identity in the previous equation, together with the fact that $|\|\bb{s}\|_1 - \|\bb{t}\|_1| \leq \|\bb{s} - \bb{t}\|_1$, one concludes that
\[
|D^{\bb{\gamma}} K_{\bb{u},v}(\bb{s}) - D^{\bb{\gamma}} K_{\bb{u},v}(\bb{t})| \leq L \, \|\bb{s} - \bb{t}\|_1^{\,\beta - m},
\]
for some constant $L \in (0, \infty)$ depending only on the integer $d$ and the reals $u_1, \ldots, u_d, v$, thereby proving the claim.

These findings are summarized below for the record; cf. Remark~2 in \cite{MR2775207}.

\begin{proposition}
\label{prop:4}
Let $K_{\bb{u},v}$ denote the $d$-variate Dirichlet density with parameters $\bb{u} = (u_1, \ldots, u_d) \in (0, \infty)^d$ and $v \in (0, \infty)$, defined in~\eqref{eq:1}.
\begin{itemize}
\item [{\rm (i)}]
If $(u_1, \ldots, u_d, v) \in \mathbb{N}^{d+1}$, then for every real $\beta \in (0, \infty)$, there exists a scalar $L \in (0, \infty)$ such that $K_{\bb{u},v} \in \Sigma (d, \beta, L)$.
\item [{\rm (ii)}]
If $(u_1, \ldots, u_d, v) \in [r,\infty)^{d+1} \setminus \mathbb{N}^{d+1}$ for some real $r \in [1, \infty)$, then for every $\beta \in (0, r - 1]$, there exists a scalar $L \in (0, \infty)$ such that $K_{\bb{u},v} \in \Sigma (d, \beta, L)$.
\end{itemize}
\end{proposition}

The somewhat broader formulation of Proposition~\ref{prop:4} owes to the fact that, in general,
\begin{equation}
\label{eq:6}
f \in \Sigma (d, \beta, L) \quad \Rightarrow \quad \forall_{\beta^\ast \in (0, \,\beta)} \;  \exists_{L^\ast \in (0, \infty)} \quad f \in \Sigma (d, \beta^\ast, L^\ast).
\end{equation}

In particular, suppose that the $L^1$ or $L^2$ loss function is preferred and that one suspects that the data at hand arise from a density $f \in \Sigma (d, \beta, L)$ with a high degree of regularity $\beta \in [2, \infty)$. Then Proposition~\ref{prop:1} could still be invoked to build an asymptotically minimax sequence of Aitchison--Lauder estimators by taking $b_n \asymp n^{-2/(d+2 \beta^\ast)}$ for some $\beta^* \in (0, 2]$ given that $f$ also belongs to the class $\Sigma (d, \beta^\ast, L^\ast)$ for some Lipschitz constant $L^\ast \in (0, \infty)$ whose value has no influence on the choice of the bandwidth parameter. Because the resulting estimator would be asymptotically minimax with respect to a much larger class of densities, however, this sleight of hand would be at the expense of the minimax rate $r_n (d, \beta, L) \asymp_{d,\beta,L} n^{-\beta /(d + 2 \beta)}$ given in~\eqref{eq:3}, which is a decreasing function of $\beta$.

\section{Conclusion\label{sec:5}}

Using results of \citet{MR4319409}, conditions were found under which the $d$-variate Dirichlet kernel density estimator defined in~\eqref{eq:2} can achieve the asymptotic minimax rate for the $L^p$ loss with a suitable choice of bandwidth parameter. As per Proposition~\ref{prop:1}, this is possible when the underlying density $f$ belongs to the H\"older space $\Sigma(d, \beta, L)$ for some Lipschitz constant $L \in (0, \infty)$ provided that the $L^p$ loss and regularity parameter $\beta$ are such that $(p,\beta)\in [1,3) \times (0,2] \cup \mathcal{A}_d$. To achieve the minimax rate, the bandwidth parameter $b \in (0, \infty)$ must vary with the sample size $n$ and depend on the degree $\beta$ of smoothness of the underlying density in such a way that $b_n \asymp n^{-2/(d + 2 \beta)}$.

It is interesting to note that because of the embedding property~\eqref{eq:6}, Silverman's rule of thumb which consists of taking $b_n \asymp n^{-2/(d+4)}$ ensures that the Aitchison--Lauder estimator is asymptotically minimax on $\Sigma (d, 2, L)$ for any $\beta$-H\"older density with regularity parameter $\beta \in [2, \infty)$. In dimension $d + 1 = 2$, for instance, this corresponds to the familiar rate $r_n (1, 2, L) \asymp_L n^{-2/5}$. Only for densities with regularity parameter $\beta \in (0, 2)$ would this be insufficient.

However, it was shown in Propositions~\ref{prop:2} and~\ref{prop:3} that the Dirichlet kernel density estimator defined in~\eqref{eq:2} cannot be asymptotically minimax on $\Sigma (d, \beta, L)$ for $L^p$ losses with $p \in [4, \infty)$ or when $\beta \in (2, \infty)$. As detailed in Remark~\ref{remark:1}, these results rectify in a minor way and, more importantly, extend to all dimensions those already reported in dimension $d + 1 = 2$ by \citet{MR2775207}. However, the case $(p, \beta) \in [3, 4) \times (0, 2] \setminus \mathcal{A}_d$ remains open.

The results reported here are generally good news for the Dirichlet kernel density estimator of \citet{doi:10.2307/2347365}. Nevertheless, there may be reasons for preferring other options. One of them is the fact that the estimator defined in~\eqref{eq:2} does not integrate to $1$, except asymptotically \citep{MR4319409}. Some users may also feel more comfortable relying on scalings of a fixed kernel function by proceeding, e.g., as proposed by \citet{doi:10.1016/j.cageo.2009.12.011}. At minima, the arguments presented here show that deriving asymptotic properties of an estimator based on a variable kernel function is not as difficult as these authors had anticipated.

\appendix

\setcounter{lemma}{0}
\setcounter{rem}{0}
\renewcommand{\thelemma}{\Alph{section}\arabic{lemma}}
\renewcommand{\therem}{\Alph{section}\arabic{rem}}

\section{Proof of Proposition~\ref{prop:1}\label{app:A}}

Fix $p \in [1, 3)$ and $\beta \in (0, 2]$. Further let $b = b_n = c \, n^{-2/(d + 2 \beta)}$ for some constant $c \in (0, \infty)$. Let $\bb{\xi}_{\bb{s}} = (\xi_1, \ldots,\xi_d)$ be a random vector having Dirichlet distribution~\eqref{eq:1} with parameters $\bb{u} = \bb{s} / b + \bb{1}$ and $v = (1 - \|\bb{s}\|_1) / b + 1$. Observe that one then has
\begin{equation}
\label{eq:A.1}
\EE \big \{ f_{n,b}(\bb{s}) \big \} = \EE \big \{ K_{\bb{s} / b + \bb{1}, (1 - \|\bb{s}\|_1) / b + 1}(\bb{X}) \big \} = \EE \big \{ f(\bb{\xi}_{\bb{s}}) \big \}.
\end{equation}
By the triangle inequality for $\| \cdot \|_p$, one also has
\[
R_n \big (f_{n,b},f \big ) = \big\{ \EE \big( \|f_{n,b} - f\|_p^p \big) \big\}^{1/p} \le \Big[ \EE \big[ \big \{ \|f_{n,b} - \EE (f_{n,b}) \|_p + \|\EE (f_{n,b}) - f\|_p \big\}^p \big] \Big]^{1/p}.
\]
Using the fact that $(t + w)^p \leq 2^{p-1} (t^p + w^p)$ for all $t, w \in [0, \infty)$, one can then deduce that
\[
R_n \big (f_{n,b},f \big ) \ll \big\{ \EE \big\{ \|f_{n,b} - \EE (f_{n,b}) \|_p^p \big\} + \|\EE (f_{n,b}) - f\|_p^p \big\}^{1/p}.
\]
Given the sub-additivity of the map $x \mapsto x^{1/p}$ on $[0, \infty)$ when $p \in [1, \infty)$, it follows that
\begin{equation}
\label{eq:A.2}
R_n \big (f_{n,b},f \big ) \ll \big[ \EE \big \{\|f_{n,b} - \EE (f_{n,b}) \|_p^p \big\} \big]^{1/p} + \big \{ \|\EE (f_{n,b}) - f\|_p^p \big\}^{1/p} \equiv A _n + B_n.
\end{equation}

\bigskip
Next, suitable bounds will be found on the terms $A_n$ and $B_n$ implicitly defined in~\eqref{eq:A.2}. The following fact will be used to bound $\sup_{\bb{s} \in \mathrm{Int}(\mathcal{S}_d)} |f(\bb{s})|$ in~\eqref{eq:A.6} and $\sup_{\bb{s} \in \mathrm{Int}(\mathcal{S}_d)} \|\nabla f(\bb{s})\|_1$ in~\eqref{eq:A.14}.

\bigskip
\begin{rem}
\emph{By generalizing the reasoning on p.~7 of the book by~\citet{MR2013911}, one can show that if $f \in \Sigma (d, \beta, L)$ for some integer $d \in \mathbb{N}$ and scalars $\beta, L  \in (0, \infty)$, then there exists a constant $M = M (d, \beta, L) \in (0, \infty)$ which is sufficiently large to ensure that
\[
\max_{\bb{\gamma} \in \mathbb{N}_0^d: \|\bb{\gamma}\|_1 \leq m} \sup_{\bb{s} \in \mathrm{Int}(\mathcal{S}_d)} | D^{\bb{\gamma}} f(\bb{s}) | \leq M.
\]}
\end{rem}

\subsection{Derivation of an upper bound on \texorpdfstring{$A_n$}{An}\label{app:A.1}}

Fix $\bb{s} \in \mathcal{S}_d$ and for each integer $i \in \{1, \ldots, n\}$, let
\begin{equation}
\label{eq:A.3}
Y_i = K_{{\bb{s}}/{b} + \bb{1}, {(1 - \|\bb{s}\|_1)}/{b} + 1}(\bb{X}_i) - \EE \big \{ f_{n,b}(\bb{s}) \big\}.
\end{equation}
Given that the observations $\bb{X}_1, \ldots, \bb{X}_n$ are mutually independent, so are the variables $Y_1, \ldots, Y_n$ and in view of~\eqref{eq:A.1}, they all have mean zero and finite variance, as does ${\bar Y}_n = (Y_1 + \cdots + Y_n)/n$.

It follows from Jensen's inequality that when $p \in [0, 2]$,
\[
\EE \big( | {\bar Y}_n |^p \big) \leq \big\{ \VV \big( {\bar Y}_n \big) \big\}^{p/2}
\]
while if $p \in (2, \infty)$, then~(4.11) of \citet{MR523165} yields
\[
\EE \big( | {\bar Y}_n |^p \big) \ll_p (\|Y_1\|_{\infty}/n)^{p-2} \VV\big( {\bar Y}_n \big) + \big\{ \VV \big( {\bar Y}_n \big) \big\}^{p/2}.
\]

In view of the definition of $Y_i$ given in~\eqref{eq:A.3}, one then has, for any $p \in [0, \infty)$ and $p^\prime = \max(0, p-2)$,
\begin{equation}
\label{eq:A.4}
\EE \big[ |f_{n,b}(\bb{s}) - \EE \{ f_{n,b}(\bb{s}) \} |^p \big] \ll_p  (\| K_{\bb{s}/b + \bb{1}, (1 - \|\bb{s}\|_1)/b + 1} \|_{\infty}/n)^{p^{\prime}} \, \VV \big \{ f_{n,b}(\bb{s}) \big\} + \big[ \VV \big\{ f_{n,b}(\bb{s}) \big\} \big]^{p/2}.
\end{equation}
Now it is already known, thanks to~(15) and Lemmas~1--2 in \cite{MR4319409} that, for $n$ large enough and $b$ small enough,
\[
\|K_{\bb{s}/b + \bb{1}, (1 - \|\bb{s}\|_1)/b + 1} \|_{\infty} \ll_d \frac{b^{-d/2} }{\sqrt{(1 - \|\bb{s}\|_1) s_1 \cdots s_d}}
\]
and
\[
\VV \big\{ f_{n,b} (\bb{s}) \big\} \ll_d \frac{n^{-1} b^{-d/2} \|f\|_{\infty} }{\sqrt{(1 - \|\bb{s}\|_1) s_1 \cdots s_d}}.
\]

By interpolating the bound on $\| K_{\bb{s}/b + \bb{1}, (1 - \|\bb{s}\|_1)/b + 1} \|_{\infty}$ above together with the uniform bound in Lemma~\ref{lemma:C3} below, one can deduce that, for any choice of real $q \in [0,1]$,
\begin{equation}
\label{eq:A.5}
\| K_{\bb{s}/b + \bb{1}, (1 - \|\bb{s}\|_1)/b + 1} \|_{\infty} \ll_d \left(\frac{b^{-d/2} }{\sqrt{(1 - \|\bb{s}\|_1) s_1 \cdots s_d}}\right)^q (b^{-d} )^{1-q}.
\end{equation}
Plugging these last two bounds into~\eqref{eq:A.4}, one concludes that
\begin{equation}
\label{eq:A.6}
\EE \big[ |f_{n,b}(\bb{s}) - \EE \{ f_{n,b}(\bb{s}) \} |^p \big]
\ll_{d, \beta, L, p} \frac{n^{-(p^\prime + 1)} (b^{-d/2})^{p^\prime q + 2 p^\prime (1 - q) + 1}}{ \{ (1 - \|\bb{s}\|_1) s_1 \cdots s_d \} ^{(p^\prime q + 1)/2}} + \frac{(n^{-1/2} \, b^{-d/4})^p}{ \{ (1 - \|\bb{s}\|_1) s_1 \cdots s_d \} ^{p/4}}.
\end{equation}

Given that the map $\bb{s} \mapsto \{ (1 - \|\bb{s}\|_1) s_1 \cdots s_d \}^{-a}$ is integrable on $\mathcal{S}_d$ if and only if $a < 1$, which imposes the restriction $p < \min(2 + 1/q, 4)$, and considering that the map $x \mapsto x^{1/p}$ is sub-additive on $[0, \infty)$ whenever $p \in [1, \infty)$, it follows from integration on both sides of~\eqref{eq:A.6} that
\[
A_n = \big[ \EE \big\{ \|f_{n,b} - \EE ( f_{n,b}) \|_p^p \big\} \big]^{1/p}
\ll_{d, \beta, L, p} \big\{ n^{-(p^\prime + 1 - p/2)/p} (b^{-d/2})^{(p^\prime q + 2 p^\prime (1 - q) + 1 - p/2)/p} + 1 \big\} \, n^{-1/2} \, b^{-d/4}.
\]
Moreover, given the assumption that $b = c n^{-2/(d + 2 \beta)}$, one has
\begin{equation}
\label{eq:A.7}
\frac{A_n}{n^{-1/2} \, b^{-d/4}} \ll_{d, \beta, L, p} n^{-(p^\prime + 1 - p/2)/p} (n^{d/(d + 2 \beta)})^{(p^\prime q + 2 p^\prime (1 - q) + 1 - p/2)/p} + 1.
\end{equation}

When $p\in [1,3)$, it suffices to take $q = 1$ for the last quantity to be bounded. When $p\in [3,4)$, one needs to be more careful because of the aforementioned restriction $p < 2 + 1/q$ or, equivalently, $q < 1 / (p - 2)$. For the right-hand side of~\eqref{eq:A.7} to be bounded, one must have
\[
-(p^\prime + 1 - p/2) + \frac{d}{d + \beta} \,  (2p^\prime - p^\prime q + 1 - p/2) < 0
\]
or equivalently
\[
- \frac{d + \beta}{d} \, (p^\prime + 1 - p/2) + (2p^\prime + 1 - p/2) < q p^\prime \quad \Leftrightarrow \quad 1 - \frac{\beta}{d p^\prime} \, (p^\prime + 1 - p/2) < q.
\]

When $p \in [3,4)$, one has $p^\prime = p - 2$, so that the above is equivalent to the simpler condition $1 - \beta/ (2d) < q$. Combined with the restriction $p < 2 + 1/q$, one then finds that $q$ must satisfy the condition
\[
1 - \frac{\beta}{2d} < q < \frac{1}{p - 2} \,
\]
and the latter can be fulfilled only if
\[
2d \, \frac{p-3}{p-2} < \beta.
\]
Given the additional restriction $\beta \leq 2$, the range of viable values of $p$ is then limited to the interval $[3, 2 + d / (d - 1)]$.

\subsection{Derivation of an upper bound on \texorpdfstring{$B_n$}{Bn}\label{app:A.2}}

It will be shown below that, for every real $\beta \in (0, 2]$, one has
\begin{equation}
\label{eq:A.8}
B_n \ll_{d, \beta, L} b^{\beta / 2}.
\end{equation}
However, the cases $\beta \in (0, 1]$ and $\beta \in (1, 2]$ must be treated separately.

\medskip
\noindent
\textit{Case $\beta \in (0, 1]$}:
For any such value of $\beta$, the map $x \mapsto x^{\beta}$ is concave. Combining this fact with identity~\eqref{eq:A.1} and the assumption that $f\in \Sigma(d, \beta, L)$, one can deduce that
\begin{equation}
\label{eq:A.9}
| \EE \{ f_{n,b}(\bb{s}) \} - f(\bb{s})| \leq \EE \big \{ |f(\bb{\xi}_{\bb{s}}) - f(\bb{s})| \big\} \leq L \, \EE \big( \|\bb{\xi}_{\bb{s}} - \bb{s}\|_1^{\beta} \big) \leq L \, \big\{ \EE \big( \|\bb{\xi}_{\bb{s}} - \bb{s}\|_1\big) \big\}^{\beta},
\end{equation}
where Jensen's inequality was used at the last step. For small enough $b$, it also follows from Jensen's inequality that
\begin{equation}
\label{eq:A.10}
\EE \big( \|\bb{\xi}_{\bb{s}} - \bb{s} \|_1 \big) = \sum_{i=1}^d \EE \big( |\xi_i - s_i| \big) \leq \sum_{i=1}^d \big\{ \EE \big( |\xi_i - s_i|^2 \big) \big\}^{1/2}.
\end{equation}

Now the fact that, for any integer $i \in \{ 1, \ldots, d\}$, the random variable $\xi_i$ has a Beta distribution with parameters $s_i/b + 1$ and $(1 - s_i)/b + d$ implies that
\begin{equation}
\label{eq:A.11}
\EE (\xi_i)  = \frac{{s_i}/{b} + 1}{{1}/{b} + d + 1}, \quad \VV (\xi_i) = \frac{({s_i}/{b} + 1) \{ {(1 - s_i)}/b + d\}}{(1/b + d + 1)^2 ({1}/{b} + d + 2)},
\end{equation}
so that for small enough real $b \in (0, \infty)$, one has
\begin{equation}
\label{eq:A.12}
\EE \big( |\xi_i - s_i|^2 \big) = \VV (\xi_i) + | \EE (\xi_i) - s_i | ^2 \leq {b}/{2} + \{b (d + 1)\} ^2 \leq b.
\end{equation}

Applying the latter bound term by term to the right-hand term of inequality~\eqref{eq:A.10}, one finds
\begin{equation}
\label{eq:A.13}
\EE \big( \|\bb{\xi}_{\bb{s}} - \bb{s} \|_1 \big) \leq d b^{1/2},
\end{equation}
and upon substitution into~\eqref{eq:A.9}, one can then conclude that $B_n \ll_{d,\beta,L} b^{\beta / 2}$, as claimed.

\bigskip
\noindent
\textit{Case $\beta \in (1, 2]$}: As $f\in \Sigma(d, \beta, L)$ by assumption, the multivariate mean value theorem implies the existence of a random vector $\bb{\zeta}_{\bb{s}}\in \mathcal{S}_d$ on the line segment joining $\bb{\xi}_{\bb{s}}$ and $\bb{s}$ such that
\begin{align}
\label{eq:A.14}
|\EE \{ f_{n, b}(\bb{s})\} - f(\bb{s})|
& \leq \EE \big\{ |f(\bb{\xi}_{\bb{s}}) - f(\bb{s})| \big\} \notag \\[2mm]
& = \EE\big[|\nabla f(\bb{s})^{\top} (\bb{\xi}_{\bb{s}} - \bb{s}) + \{ \nabla f(\bb{\zeta}_{\bb{s}}) - \nabla f(\bb{s}) \} ^{\top} (\bb{\xi}_{\bb{s}} - \bb{s})|\big] \notag \\[1mm]
&\ll_{d, \beta, L} \, \EE \big( \|\bb{\xi}_{\bb{s}} - \bb{s}\|_1 \big) + \EE \big( \|\bb{\xi}_{\bb{s}} - \bb{s}\|_1^{\beta} \big).
\end{align}

Observe that the first summand on the right-hand side of~\eqref{eq:A.14} is bounded above by $d b^{1/2}$, as shown in~\eqref{eq:A.13}. To bound the second summand from above, first apply the triangle inequality for $\|\cdot\|_1$ and then Jensen's inequality twice to get
\begin{align*}
\EE \big( \|\bb{\xi}_{\bb{s}} - \bb{s}\|_1^{\beta} \big)
& \leq \EE \big[ \big\{ \|\bb{\xi}_{\bb{s}} - \EE (\bb{\xi}_{\bb{s}}) \|_1 + \|\EE (\bb{\xi}_{\bb{s}}) - \bb{s}\|_1\big\}^{\beta} \big] \notag \\[2mm]
& \leq (2 d)^{\beta - 1} \left[ \sum_{i=1}^d \EE \big\{ |\xi_i - \EE (\xi_i)|^{\beta} \big\} + \sum_{i=1}^d |\EE (\xi_i) - s_i|^{\beta} \right] \notag \\[2mm]
& \leq (2 d)^{\beta - 1} \left [ \sum_{i=1}^d \big\{ \VV(\xi_i) \big\}^{\beta / 2} + \sum_{i=1}^d |\EE (\xi_i) - s_i|^{\beta} \right].
\end{align*}

Calling again on the fact that for each integer $i \in \{ 1, \ldots, n\}$, the random variable $\xi_i$ has a Beta distribution with parameters $s_i/b + 1$ and $(1 - s_i)/b + d$, and proceeding as in the derivation of inequality~\eqref{eq:A.12}, one finds
\begin{equation}
\label{eq:A.15}
\EE \big( \|\bb{\xi}_{\bb{s}} - \bb{s}\|_1^{\beta} \big) \ll_{d,\beta} b^{\beta / 2}.
\end{equation}
It then suffices to apply the bounds~\eqref{eq:A.13} and~\eqref{eq:A.15} to inequality~\eqref{eq:A.14} to conclude that $B_n \ll_{d, \beta, L} b^{\beta / 2}$, as claimed.

\subsection{Final step in the proof of Proposition~\ref{prop:1}\label{app:A.3}}

Upon applying bounds~\eqref{eq:A.7} and~\eqref{eq:A.8} to~\eqref{eq:A.2}, one gets
\[
R_n (f_{n,b},f) \ll A_n + B_n \ll_{d, \beta, L, p} n^{-1/2} b^{-d/4} + b^{\beta / 2} \ll_{d, \beta, L, p} \varphi_{n} (d, \beta) = n^{-\beta /(d + 2 \beta)},
\]
which, in view of~\eqref{eq:3}, yields the stated conclusion.\hfill $\Box$

\section{Proof of Propositions~\ref{prop:2}~and~\ref{prop:3}\label{app:B}}

The following functions play a role in the proofs of these results. For any $\bb{s} \in \mathcal{S}_d$, define
\[
f_0 (\bb{s}) = d!, \quad f_3 (\bb{s}) = (d - 1)! \, (d + 1) \, \|\bb{s}\|_1,
\]
and for any reals $\beta \in (0, 2]$ and $L \in (0, \infty)$ such that $0 < b \leq \min\{1 / (4d)^2, ((d! / L)^{1/\beta} / 3)^2\}$, let ${\bb{t}^{(\bb{k})}}$ be a vector whose $i$th coordinate equals $t_i^{(\bb{k})} = 1/(4d) + 3b^{1/2}(2k_i-1) $ for every integer $i \in \{ 1, \ldots, d \}$ and introduce
\begin{equation}
\label{eq:B.1}
f_{\beta}(\bb{s}) = f_{\beta,b}(\bb{s}) = d! \, \ind_{\mathcal{S}_d}(\bb{s}) + L_{\beta} \sum_{\bb{k}\in\{1,\ldots,2N\}^d} (-1)^{k_1} (3 b^{1/2})^{\beta} \psi \left(\frac{\bb{s} - \bb{t}^{(\bb{k})}}{3 b^{1/2}}\right),
\end{equation}
where $\ind_{\mathcal{S}_d}$ stands for the indicator function of the set $\mathcal{S}_d$ and
\[
L_{\beta} = L \min (1, 1/ \beta), \quad N = \lceil 1/(24 d b^{1/2}) \rceil, \quad \psi (\bb{s}) = (1-\|\bb{s}\|_1^\beta) \ind_{\mathcal{S}_d}(\bb{s}),
\]
with $\lceil x \rceil$ denoting the smallest integer greater than or equal to $x \in (0, \infty)$.

The functions $f_0$, $f_3$ and $f_{\beta} = f_{\beta,b}$ defined above are all densities with support on $\mathcal{S}_d$. This is immediate for $f_0$, considering that $\int_{\mathcal{S}_d} d! \rd \bb{s} = 1$. To prove that $f_3$ also integrates to $1$, use the change of variables $t_i = s_1 + \cdots + s_i$ for every integer $i \in \{ 1, \ldots, d \}$ to write successively
\begin{align*}
\int_{\mathcal{S}_d} \|\bb{s}\|_1 \rd \bb{s}
= \int_{0 < t_1 < \dots < t_d < 1} t_d \rd \bb{t} & = \int_{0 < t_2 < \dots < t_d < 1} t_2 t_d \rd \bb{t} \\[2mm]
& = \frac{1}{2} \int_{0 < t_3 < \dots < t_d < 1} t_3^2  t_d \rd \bb{t} = \dots = \frac{1}{(d-2)!} \int_{0 < t_{d-1} < t_d < 1} t_{d-1}^{d-2} t_d \rd \bb{t},
\end{align*}
which shows that the integral equals $1/\{ (d-1)! \, (d+1)\}$.

\begin{figure}[t!]
% FIGURE 1
\captionsetup[subfigure]{labelformat=empty}
\captionsetup{width=0.6\linewidth}
\centering

\begin{subfigure}[b]{0.4\textwidth}
\centering
\includegraphics[width=\textwidth]{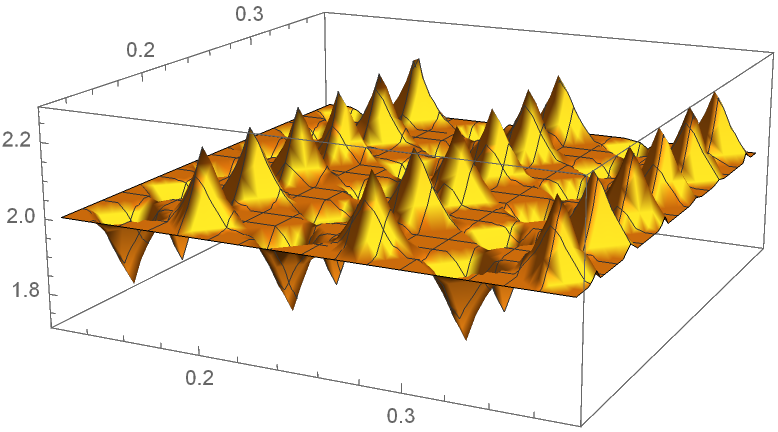}
\caption{$\beta = 1/2$ and $b = 5 \times 10^{-5}$\vspace{3mm}}
\label{fig:dirichlet.1}
\end{subfigure}
\quad
\begin{subfigure}[b]{0.4\textwidth}
\centering
\includegraphics[width=\textwidth]{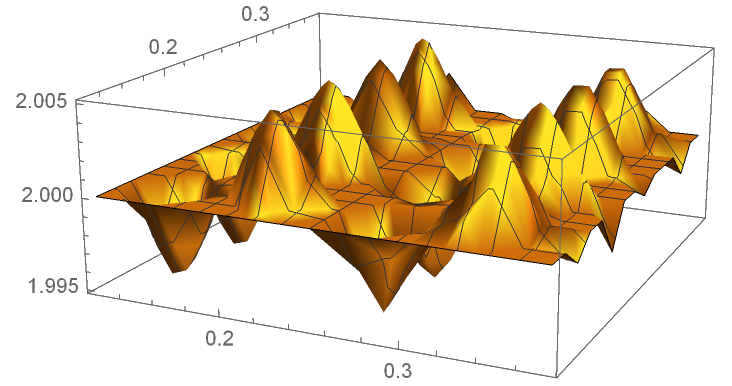}
\caption{$\beta = 3/2$ and $b = 10^{-4}$\vspace{3mm}}
\label{fig:dirichlet.2}
\end{subfigure}
\caption{Graphs of $f_{\beta}$ in the region $[1/8,3/8]^2 = [1/(4d),3/(4d)]^d$ when $d = 2$.}
\label{fig:1}
\end{figure}

Turning to $f_{\beta} = f_{\beta,b}$ defined in~\eqref{eq:B.1} for a fixed real $\beta \in (0, 2]$, note that it takes nonnegative values on $\mathcal{S}_d$, and that the volume under the positive spikes induced by the map
\[
\bb{s} \mapsto L_{\beta} (-1)^{k_1} (3 b^{1/2})^{\beta} \psi \big \{ (\bb{s} - \bb{t}^{(\bb{k})}) / (3 b^{1/2}) \big \}
\]
for even values of $k_1$ are offset by the negative spikes corresponding to odd values of $k_1$. This is illustrated in Fig.~\ref{fig:1} for two different sets of values of $(\beta, b)$ in dimension $d + 1 = 3$.

The proofs of Propositions~\ref{prop:2}~and~\ref{prop:3} rely on the following technical lemma, whose proof is deferred to~\ref{app:B.3}. In what follows, $\PP_\beta$, $\EE_\beta$ and $\VV_\beta$ respectively refer to a probability, an expectation and a variance computed with respect to density $f_\beta$, whatever $\beta \in [0, 2] \cup \{3\}$.

\begin{lemma}
\label{lemma:B1}
The following statements hold true for every real $b \in (0,1)$.
\begin{enumerate}[(a)]
\item
For all reals $p \in [2, \infty)$ and $L \in (0, \infty)$, one has $f_0 \in \Sigma(d, \beta, L)$ and
\[
\EE_0 \big(\|f_{n,b} - f_0\|_p^p \big) \gg_{d,p} \, \frac{\mathcal{I}_n(b,p)}{(n \, b^{d/2})^{p/2}},
\]
where
\[
\mathcal{I}_n(b,p) =  \int_{\mathcal{S}_d(b)} \{ (1 - \|\bb{s}\|_1) s_1 \cdots s_d \} ^{-p/4} \rd \bb{s} \quad and \quad
\mathcal{S}_d(b) =  \big\{\bb{s}\in \mathcal{S}_d: 1 - \|\bb{s}\|_1 \geq b ~\mbox{\rm and}~ \min (s_1, \ldots, s_d) \geq b \big\}.
\]

\item
For all reals $p \in [2, \infty)$, $\beta \in (0,2]$, and $L \in (0, \infty)$, one has $f_{\beta} \in \Sigma(d, \beta, L)$ and
\[
\EE_{\beta} \big( \|f_{n,b}-f_{\beta}\|_p^p \big) \gg_{d, \beta, p} b^{p \beta / 2},
\]
\item
For all reals $p \in [1, \infty)$, $\beta \in (2, \infty)$, and $L \in (0, \infty)$, one has $f_3 \in \Sigma(d, \beta, L)$ and
\[
\EE_3 \big( \|f_{n,b}-f_3\|_p^p \big) \gg_{d,p} b^{p}.
\]
\end{enumerate}
\end{lemma}

\subsection{Proof of Proposition~\ref{prop:2} assuming Lemma~\ref{lemma:B1}\label{app:B.1}}

Fix reals $p \in [1, \infty)$, $\beta \in (2, \infty)$, and $L \in (0, \infty)$. To establish the result, first observe that
\[
R_n(f_{n,b}, f_3) \ge \EE_3 \big( \|f_{n,b} - f_3\|_1\big).
\]
Indeed, using Fubini's theorem and Jensen's inequality twice, one finds
\begin{align*}
R_n(f_{n,b}, f_3) = \left[ \int_{\mathcal{S}_d} \EE_3 \big\{ |f_{n,b}(\bb{s}) - f_3(\bb{s})|^p \big\} \rd \bb{s}\right]^{1/p}
& \geq (d!)^{1 - 1/p} \int_{\mathcal{S}_d} \left[ \EE_3 \big\{ |f_{n,b}(\bb{s}) - f_3(\bb{s})|^p \big\} \right]^{1/p} \rd \bb{s} \\[2mm]
& \geq \int_{\mathcal{S}_d} \EE_3 \big\{ |f_{n,b}(\bb{s}) - f_3(\bb{s})| \big\} \rd \bb{s} = \EE_3 \big( \|f_{n,b} - f_3\|_1\big).
\end{align*}
Therefore, to conclude the proof, it is sufficient to show that
\begin{equation}
\label{eq:B.2}
\liminf_{n\to \infty} \frac{\EE_3 \big( \|f_{n,b} - f_3\|_1 \big)} {\varphi_n (d, \beta)} = \infty,
\end{equation}
where $\varphi_n (d, \beta)$ is defined as in~\eqref{eq:3}. If $b = b_n$ is such that
\[
\liminf_{n\to \infty} \EE_3 (\|f_{n,b} - f_3\|_1) > 0,
\]
then the conclusion of Proposition~\ref{prop:2} follows trivially. Therefore, assume for the remainder of the proof that
\[
\liminf_{n\to \infty} \EE_3 (\|f_{n,b} - f_3\|_1) = 0
\]
and let $b = b_{k_n}$ be any subsequence that achieves this $\liminf$.

By the asymmetric kernel analog of Theorem~1~[(i) $\Rightarrow$ (v)] in \cite{MR707939}, with the identification of bandwidths $b = h^2$, one must have $b\to 0$ and $n b^{d/2}\to \infty$ as $n\to \infty$, which makes the relation $b = \oo_d(n^{-2/d})$ impossible.

With this in mind, one can divide the proof of~\eqref{eq:B.2} into the two remaining cases, namely
\begin{itemize}
\item [(i)]
$b \gg_d n^{-2/(d + 4)}$;
\item [(ii)]
$n^{-2/d} \ll_d b = \oo_d(n^{-2/(d + 4)})$.
\end{itemize}

\noindent
\textit{Case $(i)$}: By part (c) of Lemma~\ref{lemma:B1} with $p = 1$, one has
\[
\EE_3 (\|f_{n,b} - f_3\|_1) \gg_d b \gg_d n^{-2/(d + 4)} = \varphi_n (d, \beta) {\varphi_n (d, 2)}/{\varphi_n (d, \beta)}.
\]
Given that $\lim_{n\to \infty} {\varphi_n (d, 2)}/{\varphi_n (d, \beta)} = \infty$ when $\beta > 2$, the argument is complete.

\bigskip
\noindent
\textit{Case $(ii)$}: The proof here follows the same argument as in the second case considered in the proof of Theorem~1 in \cite{MR760686}; see pp.~1259--1260 in that paper. First note that by the triangle inequality and Fubini's theorem, one has
\[
\EE_3 \big(\|f_{n,b} - f_3\|_1 \big)
\geq \frac{1}{2} \, \EE_3 \left[ \int_{\mathcal{S}_d} |f_{n,b}(\bb{s}) - \EE_3 \{ f_{n,b}(\bb{s}) \} | \rd \bb{s }\right] = \frac{1}{2} \int_{\mathcal{S}_d} \EE_3 \left[ |f_{n,b}(\bb{s}) - \EE_3 \{ f_{n,b}(\bb{s}) \} |\right] \rd \bb{s}.
\]
Next, invoke Fatou's lemma to deduce that
\begin{equation}
\label{eq:B.3}
\liminf_{n\to \infty} n^{2/(d + 4)} \, \EE_3 \big(\|f_{n,b} - f_3\|_1\big) \geq \frac{1}{2} \int_{\mathcal{S}_d} \liminf_{n\to \infty} n^{2/(d + 4)} \, \EE_3 \left[ |f_{n,b}(\bb{s}) - \EE_3 \{ f_{n,b}(\bb{s}) \} |\right] \rd \bb{s}.
\end{equation}

Now, by Markov's inequality, one has, for any real $C \in (0, \infty)$,
\[
n^{2/(d + 4)} \, \EE_3 [|f_{n,b}(\bb{s}) - \EE_3 \{ f_{n,b}(\bb{s}) \} | ] \geq C \, \PP_3 [ |f_{n,b}(\bb{s}) - \EE_3 \{ f_{n,b}(\bb{s}) \} | \geq \frac{C}{n^{2/(d + 4)}} ].
\]
Also, by the Berry--Esseen theorem and the fact that, for all $u, v, w \in \mathbb{R}^3$, $|u - v| \leq w$ implies $u \geq v - w$, the right-hand term of the above inequality is bounded below by
\[
C \, \PP_3 \left[ |Z| \geq \frac{C}{n^{2/(d + 4)} \sqrt{\VV_3 \{ f_{n,b}(\bb{s})\} }}\right] - \OO_d \left[ \frac{\EE_3 \{ |\zeta_1 - \EE_3 (\zeta_1)|^3 \} }{n^{1/2} \{\VV_3(\zeta_1)\}^{3/2}}\right],
\]
where $Z$ stands for a standard $\mathcal{N}(0,1)$ Gaussian random variable and $\zeta_1 = K_{\bb{s} / b + \bb{1}, (1 - \|\bb{s}\|_1) / b + 1}(\bb{X}_1)$.

It was already shown by \citet{MR4319409} in the proof of their Theorem~3 that, for every vector $\bb{s}\in \mathrm{Int}(\mathcal{S}_d)$,
\[
\frac{ \EE_3 \{ |\zeta_1 - \EE_3 (\zeta_1) |^3\} } { \VV_3(\zeta_1)} = \OO_{d,\bb{s}}(b^{-d/2}),
\]
and in their Theorem~1 that, for every such vector $\bb{s} \in \mathrm{Int}(\mathcal{S}_d)$,
\[
\VV_3 \{ f_{n,b}(\bb{s})\} = n^{-1} \VV_3(\zeta_1) = n^{-1} b^{-d/2} \{\psi(\bb{s}) f_3(\bb{s}) + \OO_{d,\bb{s}}(b^{1/2}) \}.
\]

Assuming that $b = \oo_d(n^{-2/(d + 4)})$, one can then deduce that, for every vector $\bb{s}\in \mathrm{Int}(\mathcal{S}_d)$,
\begin{align*}
n^{2/(d + 4)} \, \EE_3 [  | f_{n,b}(\bb{s}) - \EE_3 \{ f_{n,b}(\bb{s}) \}  |  ]
& \geq C \, \PP_3 \{ |Z| \geq C \, \oo_{d,\bb{s}}(1) \} - \OO_{d,\bb{s}}(n^{-1/2} b^{-d/4}) \\[2mm]
& \geq C \{ 1 + \oo_{d,\bb{s}}(1)\} - \OO_{d,\bb{s}}(n^{-1/2} b^{-d/4}).
\end{align*}
Under the assumption that $n^{-2/d} \ll_d b$ and by letting $C\to \infty$, one then has, for every vector $\bb{s}\in \mathrm{Int}(\mathcal{S}_d)$,
\[
\liminf_{n\to \infty} n^{2/(d + 4)} \, \EE_3 [ | f_{n,b}(\bb{s}) - \EE_3 \{ f_{n,b}(\bb{s}) \} | ] = \infty.
\]
Therefore, one can conclude from~\eqref{eq:B.3} that
\[
\EE_3 \big( \|f_{n,b} - f_3\|_1\big) \gg_d n^{-2/(d + 4)} = \varphi_n (d, \beta) \varphi_n(d,2) / \varphi_n (d, \beta).
\]
Given that $\lim_{n\to \infty} {\varphi_n (d, 2)}/{\varphi_n (d, \beta)} = \infty$ when $\beta \in (2, \infty)$, the argument is complete.

\bigskip
Cases (i) and (ii) having been successfully dealt with, Proposition~\ref{prop:2} is established.\hfill $\Box$

\subsection{Proof of Proposition~\ref{prop:3} assuming Lemma~\ref{lemma:B1}\label{app:B.2}}

Fix reals $p \in [4, \infty)$, $\beta \in (0, 2]$, and $L \in (0, \infty)$. Using parts (a) and (b) of Lemma~\ref{lemma:B1}, one has
\begin{align*}
\sup_{f \in \Sigma(d, \beta, L)} R_n(f_{n,b}, f)
& \geq \frac{1}{2} \big\{ R_n(f_{n,b}, f_0) + R_n(f_{n,b}, f_{\beta}) \big\} \\[1mm]
& \gg_{d, \beta, p} \inf_{b\in (0,1)} \left[ \frac{\{\mathcal{I}_n(b,p) \}^{1/p}}{(n b^{d/2})^{1/2}} + b^{\beta / 2}\right]
 \gg_{d, \beta, p} \inf_{b\in (0,1)} \left\{ \frac{|\ln b|^{1/p}}{(n b^{d/2})^{1/2}} + b^{\beta / 2}\right\}.
\end{align*}

Next, simple calculus can be used to check that whenever $n$ is larger than a certain threshold, the map
\[
b \mapsto |\ln b|^{1/p} /(n b^{d/2})^{1/2} + b^{\beta / 2}
\]
attains its minimum $b_n^0$ on $(0,1)$ and that this minimum is such that
\[
b_n^0 = c_1 n^{-2/(d + 2\beta)} \{ c_2 |\ln b_n^0|^{1/p - 1} + c_3 |\ln b_n^0|^{1/p}\}^{4/(d + 2\beta)},
\]
where $c_1 = (2/\beta)^{4/(d + 2\beta)}$, $c_2 = 1/p$ and $c_3 = d/4$.

Moreover, it can be proved easily that $b_n^0\to 0$ as $n\to \infty$. Consequently,
\[
\inf_{b\in (0,1)} \left\{\frac{|\ln b|^{1/p}}{(n b^{d/2})^{1/2}} + b^{\beta / 2}\right\} = \frac{|\ln b_n^0|^{1/p}}{\{ n (b_n^0)^{d/2}\}^{1/2}} + (b_n^0)^{\beta / 2} \gg_{d, \beta, p} \varphi_n (d, \beta) \times |\ln (b_n^0)|^{{2 \beta}/\{p (d + 2 \beta)\}}.
\]
The conclusion of Proposition~\ref{prop:3} then follows at once.\hfill $\Box$

\subsection{Proof of Lemma~\ref{lemma:B1}\label{app:B.3}}

\noindent
\textbf{Proof of Part (a)}: After trivial adjustments to the proof of Lemma~6 in \cite{MR2775207}, one can assert that
\begin{equation}
\label{eq:B.4}
\EE_0 \big( \|f_{n,b}-f_0\|_p^p \big) \geq 2^{-p} \int_{\mathcal{S}_d} \big[ \VV_0 \{ f_{n,b}(\bb{s})\} \big] ^{p/2} \rd \bb{s}.
\end{equation}
Now, \citet{MR4319409} (see p.\,13 therein) showed that, for any target density $f$ whose support is on $\mathcal{S}_d$, one has, for all $\bb{s} \in \mathcal{S}_d$,
\begin{equation}
\label{eq:B.5}
\VV_0 \{ f_{n,b}(\bb{s})\} = n^{-1} A_b(\bb{s}) \, \EE\{f(\bb{\gamma}_{\bb{s}})\} - \OO(n^{-1}),
\end{equation}
where $\bb{\gamma}_{\bb{s}}$ is random vector with distribution $\mathrm{Dirichlet}\hspace{0.3mm}(2 \bb{s} / b + \bb{1}, 2 (1 - \|\bb{s}\|_1) / b + 1)$, and
\begin{multline}
\label{eq:B.6}
A_b(\bb{s}) = \frac{b^{(d + 1) / 2} \, (1 / b + d)^{d + 1/2}}{(4\pi)^{d/2} \sqrt{(1 - \|\bb{s}\|_1) s_1 \cdots s_d}} \, \bigg(\frac{2 / b + 2d}{2 / b + d}\bigg)^{2/b + d + 1/2} e^{-d}
\times \frac{R^2 \{ (1 - \|\bb{s}\|_1) / b\} \prod_{i=1}^d R^2(s_i / b)}{R \{ 2(1 - \|\bb{s}\|_1) / b\} \prod_{i=1}^d R(2 s_i / b)} \, \frac{R(2 / b + d)}{R^2(1 / b + d)},
\end{multline}
where, for any real $z \in [0, \infty)$,
\begin{equation}
\label{eq:B.7}
R(z) = \frac{\sqrt{2\pi} }{\Gamma(z + 1)}\, e^{-z} z^{z + 1/2}.
\end{equation}

As is well-known, the map $z \mapsto R(z)$ is increasing on $(1, \infty)$, and $R(z) < 1$ for every real $z \in [1, \infty)$; see, e.g., Theorem~2.2 of~\citet{MR3684463}. Accordingly, one has, for any vector $\bb{s} \in \mathcal{S}_d$,
\[
\frac{R^2 \{ (1 - \|\bb{s}\|_1) / b\} \prod_{i=1}^d R^2(s_i / b)}{R \{ 2(1 - \|\bb{s}\|_1) / b\} \prod_{i=1}^d R(2 s_i / b)} \, \frac{R(2 / b + d)}{R^2(1 / b + d)} \geq R^{2(d+1)}(1) \, \ind_{\mathcal{S}_d(b)}(\bb{s}).
\]
By plugging the above bound into~\eqref{eq:B.6}, one deduces that, for every vector $\bb{s} \in \mathcal{S}_d$,
\begin{equation}
\label{eq:B.8}
A_b(\bb{s}) \gg_d \frac{b^{-d/2}}{\sqrt{(1 - \|\bb{s}\|_1) s_1 \cdots s_d}} \, \ind_{\mathcal{S}_d(b)}(\bb{s}).
\end{equation}
In view of~\eqref{eq:B.8}, one can also deduce from~\eqref{eq:B.5} with $f = f_0$ that, for every vector $\bb{s} \in \mathcal{S}_d(b)$,
\[
\VV_0 \big \{ f_{n,b}(\bb{s}) \big\} \gg_d \frac{n^{-1} b^{-d/2}}{\sqrt{(1 - \|\bb{s}\|_1) s_1 \cdots s_d}}.
\]
The conclusion then follows from~\eqref{eq:B.4}.\hfill $\Box$

\bigskip
\noindent
\textbf{Proof of Part (b)}: For any vector $\bb{k} \in \{1, \ldots, 2N\}^d$ and reals $b\in (0,1)$, $\e \in (0, 1/2]$, define the sets
\[
T_{\bb{k}}(\e,b) = \{\bb{t}\in \mathbb{R}^d : \|\bb{t} - \bb{t}^{(\bb{k})}\|_1 \leq \e b^{1/2}\}, \quad I_{\bb{k}}(b) = \{\bb{t}\in \mathbb{R}^d : b^{1/2} \leq \|\bb{t} - \bb{t}^{(\bb{k})}\|_1 \leq 2 b^{1/2}\}.
\]
Note that for any real $p \in [1, \infty)$, it follows from Fubini's theorem and Jensen's inequality that
\begin{align*}
\EE_{\beta} \big( \|f_{n,b}-f_{\beta}\|_p^p\big) & = \int_{\mathcal{S}_d} \EE_{\beta} \big\{ |f_{n,b}(\bb{s}) - f_{\beta}(\bb{s})|^p \big\} \rd \bb{s}
\geq \int_{\mathcal{S}_d} \big|\EE_{\beta} \{ f_{n,b}(\bb{s})\} - f_{\beta}(\bb{s})\big|^p \rd \bb{s} \notag \\
& \geq \left[ \int_{\mathcal{S}_d} \big|\EE_{\beta} \{ f_{n,b}(\bb{s})\} - f_{\beta}(\bb{s})\big| \rd \bb{s}\right]^p.
\end{align*}

Let $\Delta_N = \big \{ \bb{k}\in\{1,\ldots,2N\}^d :k_1 = 2\ell_1 ~\text{with}~ \ell_1\in\{1,\ldots,N\} \big \}$. Then
\begin{align}
\label{eq:B.9}
\int_{\mathcal{S}_d} \big|\EE_{\beta} \{ f_{n,b}(\bb{s}) \} - f_{\beta}(\bb{s})\big| \rd \bb{s}
& \geq \sum_{\bb{k}\in \Delta_N } \int_{T_{\bb{k}}(\e,b)} \big| \EE_{\beta} \{ f_{n,b}(\bb{s})\} - f_{\beta}(\bb{s})\big| \rd \bb{s} \notag \\
& \geq \sum_{\bb{k}\in \Delta_N}  \int_{T_{\bb{k}}(\e,b)}\left|\int_{\mathcal{S}_d} K_{\bb{s}/b + \bb{1}, {(1 - \|\bb{s}\|_1)}/{b} + 1}(\bb{u}) \{ f_{\beta}(\bb{u})-f_{\beta}(\bb{s}) \}\rd \bb{u}\right| \rd \bb{s} \notag \\
& \geq \sum_{\bb{k}\in \Delta_N}  \int_{T_{\bb{k}}(\e,b)} \left\{ A_{\bb{k}}(\bb{s})- B(\bb{s})\right\} \rd \bb{s},
\end{align}
where, for every vector $\bb{s} \in \mathcal{S}_d$,
\[
A_{\bb{k}} (\bb{s}) = \int_{I_{\bb{k}}(b)} K_{\bb{s}/b + \bb{1}, {(1 - \|\bb{s}\|_1)}/{b} + 1}(\bb{u}) \{ f_{\beta}(\bb{s}) - f_{\beta}(\bb{u})\} \rd \bb{u},
\]
and
\[
B(\bb{s}) = \int_{f_{\beta}(\bb{u})\geq f_{\beta}(\bb{s})} K_{\bb{s}/b + \bb{1}, {(1 - \|\bb{s}\|_1)}/{b} + 1}(\bb{u}) \{ f_{\beta}(\bb{u}) - f_{\beta}(\bb{s})\} \rd \bb{u}.
\]
Now for arbitrary vectors $\bb{k}\in \Delta_N$, $\bb{s}\in T_{\bb{k}}(\e,b)$ and $\bb{u}\in I_{\bb{k}}(b)$, one has
\[
f_{\beta}(\bb{s}) - f_{\beta}(\bb{u}) = L_{\beta}(\|\bb{u} - \bb{t}^{(\bb{k})}\|_1^{\beta} - \|\bb{s} - \bb{t}^{(\bb{k})}\|_1^{\beta}) \geq L_{\beta} b^{\beta/2} (1-\e^{\beta/2}) \geq L_{\beta} b^{\beta/2} \{ 1-(1/2)^{\beta/2} \}.
\]
Therefore,
\[
A_{\bb{k}}(\bb{s}) \geq L_{\beta} \{ 1-(1/2)^{\beta/2} \} b^{\beta/2}\int_{I_{\bb{k}}(b)} K_{\bb{s}/b + \bb{1}, {(1 - \|\bb{s}\|_1)}/{b} + 1}(\bb{u}) \rd \bb{u}.
\]

\bigskip
Next, one must call on the following lemma, whose proof is deferred to~\ref{app:B.4}.

\bigskip
\begin{lemma}
\label{lemma:B2}
Fix $\bb{s}\in [1/(4d),3/(4d)]^d$ and $\delta \in (0, 3)$. Then, as $b \to 0$, $K_{\bb{s}/b + \bb{1}, {(1 - \|\bb{s}\|_1)}/{b} + 1}(\bb{s} + \delta b^{1/2}) \gg_d b^{-d/2}$.
\end{lemma}

\bigskip
Note that for any $\bb{s}\in T_{\bb{k}}(\e,b)$ and $\bb{u}\in I_{\bb{k}}(b)$, one has $\|\bb{s} - \bb{u}\|_1 < 3b^{1/2}$. Using Lemma~\ref{lemma:B2}, one can then deduce that
\[
A_{\bb{k}}(\bb{s}) \gg_d b^{\beta/2} b^{-d/2}\int_{I_{\bb{k}}(b)} \rd \bb{u} \gg_d b^{\beta/2} b^{-d/2} b^{d/2} = b^{\beta/2}.
\]
Hence for a strictly positive constant $c_1$ that only depends on $d$, one has
\begin{equation}
\label{eq:B.10}
A_{\bb{k}}(\bb{s})\geq c_1 b^{\beta/2}.
\end{equation}

Now, considering that $f_{\beta}(\bb{s}) > d!$ whenever $\bb{s} \in T_{\bb{k}}(\e,b)$ for some vector $\bb{k} \in \Delta_N$, one has
\begin{align}\label{eq:B.11}
B (\bb{s})
&= L_{\beta} \int_{f_{\beta}(\bb{u})\geq f_{\beta}(\bb{s})} K_{\bb{s}/b + \bb{1}, {(1 - \|\bb{s}\|_1)}/{b} + 1}(\bb{u}) \Bigg[ \|\bb{s}-\bb{t}^{(\bb{k})}\|_1^\beta- \sum_{\bb{\ell}\in\Delta_N} \|\bb{u}-\bb{t}^{(\bb{\ell})}\|_1^\beta \ind_{\mathcal{S}_d}\bigg\{\frac{\bb{u}-\bb{t}^{(\bb{\ell})}}{3 b^{1/2}}\bigg\}\Bigg] \rd \bb{u} \notag \\[2mm]
& \leq L_{\beta} \e^\beta b^{\beta/2}\int_{\mathcal{S}_d} K_{\bb{s}/b + \bb{1}, {(1 - \|\bb{s}\|_1)}/{b} + 1}(\bb{u}) \rd \bb{u} = L_{\beta} \e^\beta b^{\beta/2}.
\end{align}

Finally, applying the bounds given in~\eqref{eq:B.10}--\eqref{eq:B.11} to~\eqref{eq:B.9}, one finds that
\begin{align*}
\sum_{\bb{k}\in \Delta_N}  \int_{T_{\bb{k}}(\e,b)} \{ A_{\bb{k}}(\bb{s})- B(\bb{s}) \} \rd \bb{s}
& \geq (c_1 b^{\beta/2}-L_{\beta} \e^\beta b^{\beta/2} ) \frac{1}{d!} \sum_{\bb{k}\in \Delta_N} (\e b^{1/2})^d \\[2mm]
& \geq (c_1 b^{\beta/2}-L_{\beta} \e^\beta b^{\beta/2} ) \left(\frac{b^{-1/2}}{24d}\right)^d \frac{1}{d!} (\e b^{1/2})^d.
\end{align*}
Taking $\e = \min \big [ 1/2, (c_1 / \{ 2L_{\beta}) \}^{1/\beta} \big]$, one can then deduce that
\[
\sum_{\bb{k}\in \Delta_N}  \int_{T_{\bb{k}}(\e,b)} \{ A_{\bb{k}}(\bb{s})- B(\bb{s}) \} \rd \bb{s} \gg_{d,\beta} b^{\beta/2},
\]
which leads to the statement of Part (b) via~\eqref{eq:B.9}.\hfill $\Box$

\bigskip
\noindent
\textbf{Proof of Part (c)}: Fubini's theorem and Jensen's inequality entail that, for any real $p \in [1, \infty)$,
\begin{equation}
\label{eq:B.12}
\EE_3 \big(\|f_{n,b} - f_3\|_p^p \big) = \int_{\mathcal{S}_d} \EE_3 \big\{ |f_{n,b}(\bb{s}) - f_3(\bb{s})|^p \big\} \rd \bb{s} \geq \int_{\mathcal{S}_d} |\EE_3 \{ f_{n,b}(\bb{s}) \} - f_3(\bb{s}) |^p \rd \bb{s}.
\end{equation}

To find a lower bound on the right-hand term, consider $\bb{s} \in \mathcal{S}_d$ and let $\bb{\xi}_{\bb{s}} = (\xi_1, \dots,\xi_d)$ be a random vector with distribution $\mathrm{Dirichlet}\hspace{0.3mm}(\bb{s} / b + \bb{1}, (1 - \|\bb{s}\|_1) / b + 1)$. It follows from the expectation estimate given in~\eqref{eq:A.11} that
\begin{equation}
\label{eq:B.13}
\EE_3 \{ f_{n,b}(\bb{s})\} - f_3(\bb{s}) = (d-1)! \, (d+1) \sum_{i=1}^d \big\{ \EE (\xi_i) - s_i \big\} = (d-1)! \, (d+1) \, \frac{b \{d - \|\bb{s}\|_1 (d + 1)\}}{1 + b(d+1)} .
\end{equation}
By applying identity~\eqref{eq:B.13} to~\eqref{eq:B.12}, one then finds
\[
\EE_3 \big(\|f_{n,b}-f_3\|_p^p\big) \gg_{d,p} b^{p}.
\]
This concludes the proof of Part (c) and of Lemma~\ref{lemma:B1}.\hfill $\Box$

\subsection{Proof of Lemma~\ref{lemma:B2}\label{app:B.4}}

For arbitrary vector $\bb{s} \in \mathcal{S}_d$ and reals $b, \delta \in (0, \infty)$, write
\[
K_{{\bb{s}}/{b} + \bb{1}, {(1 - \|\bb{s}\|_1)}/{b} + 1}(\bb{s} + \delta b^{1/2})
= \frac{\Gamma(1/b + d + 1)}{\Gamma \{ (1 - \|\bb{s}\|_1)/b + 1\} \prod_{i=1}^d \Gamma(s_i/b + 1)} \, (1 - \|\bb{s} + \delta b^{1/2}\|_1)^{(1 - \|\bb{s}\|_1)/b} \prod_{i=1}^d (s_i + \delta b^{1/2})^{s_i/b}.
\]
To find a lower bound on this expression, first note that it can be written in the form $K_{{\bb{s}}/{b} + \bb{1}, {(1 - \|\bb{s}\|_1)}/{b} + 1}(\bb{s}) Q_{b, \delta} (\bb{s})$, where
\[
Q_{b, \delta} (\bb{s}) = \left \{ 1 - \frac {d \delta b^{1/2}} {(1 - \|\bb{s}\|_1)} \right\}^{(1 - \|\bb{s}\|_1)/b} \prod_{i=1}^d \big( 1 + {\delta b^{1/2}}/{s_i} \big)^{s_i/b},
\]
and that one has, uniformly for $\bb{s} \in [1/(4d), 3/(4d)]^d$,
\begin{align*}
Q_{b, \delta} (\bb{s}) & = \exp \left[ \sum_{i=1}^d \frac{s_i}{b} \left( \frac{\delta b^{1/2}}{s_i}-\frac{\delta^2 b}{2s_i^2} \right) + \frac{(1 - \|\bb{s}\|_1)}{b} \left\{ -\frac{d\delta b^{1/2}}{(1 - \|\bb{s}\|_1)}-\frac{d^2\delta^2 b}{2(1 - \|\bb{s}\|_1)^2} \right\} + \OO(b^{1/2}) \right] \\[2mm]
& = \exp\left\{-\frac{\delta^2}{2}\sum_{i=1}^d\frac{1}{s_i} -\frac{d^2\delta^2}{2(1 - \|\bb{s}\|_1)}+\OO(b^{1/2})\right\} \gg_d 1.
\end{align*}
Next, for arbitrary vector $\bb{s} \in \mathcal{S}_d$ and reals $b, \delta \in (0, \infty)$, set
\[
R_b (\bb{s}) = \frac{R \{ (1 - \|\bb{s}\|_1)/b + 1\}}{R(1/b + d + 1)} \prod_{i=1}^d R(s_i/b + 1),
\]
where the function $R$ is as defined in~\eqref{eq:B.7}. It can then be checked by substitution that, for every vector $\bb{s} \in \mathcal{S}_d$ and real $b \in (0, \infty)$,
\[
K_{{\bb{s}}/{b} + \bb{1}, {(1 - \|\bb{s}\|_1)}/{b} + 1}(\bb{s}) = R_b (\bb{s}) W_b (\bb{s}),
\]
where
\begin{align*}
W_b (\bb{s}) &= \frac{\sqrt{2\pi} e^{-(1/b + d)} (1/b + d)^{1/b + d + 1/2}}{\sqrt{2\pi} e^{-(1 - \|\bb{s}\|_1)/b} \{ (1 - \|\bb{s}\|_1)/b\}^{(1 - \|\bb{s}\|_1)/b + 1/2} \prod_{i=1}^d \sqrt{2\pi} e^{-s_i/b} (s_i/b)^{s_i/b + 1/2}} \, (1 - \|\bb{s}\|_1)^{(1 - \|\bb{s}\|_1)/b} \prod_{i=1}^d s_i^{s_i/b} \\[3mm]
& = \frac{e^{-d} (1/b + d)^{1/b + d + 1/2} b^{(1 - \|\bb{s}\|_1)/b + 1/2} \prod_{i=1}^d b^{s_i/b + 1/2}}{(2\pi)^{d/2} \sqrt{(1 - \|\bb{s}\|_1) \prod_{i=1}^d s_i}} = \frac{b^{-d/2} e^{-d} (1 + bd)^{1/b + d + 1/2}}{(2\pi)^{d/2} \sqrt{(1 - \|\bb{s}\|_1) \prod_{i=1}^d s_i}}.
\end{align*}
Consequently,
\[
W_b (\bb{s})\sim_{d,\bb{s}} \frac{b^{-d/2}}{(2\pi)^{d/2} \sqrt{(1 - \|\bb{s}\|_1) \prod_{i=1}^d s_i}}
\]
and hence, as $b \to 0$,
\[
K_{{\bb{s}}/{b} + \bb{1}, {(1 - \|\bb{s}\|_1)}/{b} + 1}(\bb{s}) \gg b^{-d/2}
\]
uniformly for $\bb{s}\in[1/(4d), 3/(4d)]^d$. This completes the argument.\hfill $\Box$

\section{A uniform bound on the Dirichlet kernel}

The lemma below gives an upper bound on the Dirichlet kernel $\bb{x}\mapsto K_{{\bb{s}}/{b} + \bb{1}, {(1 - \|\bb{s}\|_1)}/{b} + 1}(\bb{x})$ from~\eqref{eq:2} which is uniform in $\bb{s}$ and $\bb{x}$ on the simplex $\mathcal{S}_d$. This result generalizes to all dimensions the analogous result for the Beta kernel ($d = 1$) stated in~(A.11) of \citet{MR1742101}.

\begin{lemma}\label{lemma:C3}
For every integer $d\in \N$, every real $b \in (0, \infty)$, and all vectors $\bb{s}\in \mathcal{S}_d$, one has
\[
\|K_{{\bb{s}}/{b} + \bb{1}, {(1 - \|\bb{s}\|_1)}/{b} + 1}\|_{\infty} = \max_{\bb{x}\in \mathcal{S}_d} K_{{\bb{s}}/{b} + \bb{1}, {(1 - \|\bb{s}\|_1)}/{b} + 1}(\bb{x}) \leq \prod_{i=1}^d (1 / b + i).
\]
\end{lemma}

\begin{proof}[\bf Proof]
First, it is well known that for the Dirichlet density $\bb{x}\mapsto K_{\bb{u},v}(\bb{x})$ as defined in~\eqref{eq:1}, the mode is attained at $\bb{x} = (\bb{u} - \bb{1}) / (\|\bb{u}\|_1 + v - d - 1)$; see, e.g., Theorem~2.4 of \citet{MR2830563}. Of interest here is the case where $\bb{u} = \bb{s} / b + \bb{1}$ and $v = (1 - \|\bb{s}\|_1) / b + 1$, so that the mode of the kernel $\bb{x}\mapsto K_{{\bb{s}}/{b} + \bb{1}, {(1 - \|\bb{s}\|_1)}/{b} + 1}(\bb{x})$ is attained exactly at $\bb{x} = \bb{s}$. Then
\begin{align*}
\max_{\bb{x}\in \mathcal{S}_d} K_{{\bb{s}}/{b} + \bb{1}, {(1 - \|\bb{s}\|_1)}/{b} + 1}(\bb{x})
&= K_{{\bb{s}}/{b} + \bb{1}, {(1 - \|\bb{s}\|_1)}/{b} + 1}(\bb{s}) \\[-1mm]
&= \frac{\Gamma(1 / b + d + 1)}{\Gamma \{(1 - \|\bb{s}\|_1) / b + 1\} \prod_{i=1}^d \Gamma(s_i / b + 1)} \, (1 - \|\bb{s}\|_1)^{(1 - \|\bb{s}\|_1) / b} \prod_{i=1}^d s_i^{s_i / b}.
\end{align*}

The goal is to maximize this last expression in $\bb{s}$. A nice feature of the function $\bb{s}\mapsto K_{{\bb{s}}/{b} + \bb{1}, {(1 - \|\bb{s}\|_1)}/{b} + 1}(\bb{s})$ is that it is log-convex on the simplex $\mathcal{S}_d$. Indeed, upon taking the logarithm and differentiating twice, one finds that, for every integers $i, j \in \{1, \ldots, d \}$ and all vectors $\bb{s}\in \mbox{Int} (\mathcal{S}_d)$,
\[
\frac{\partial^2}{\partial s_i \partial s_j} \ln \{ K_{{\bb{s}}/{b} + \bb{1}, {(1 - \|\bb{s}\|_1)}/{b} + 1}(\bb{s}) \} = \frac{1}{b^2} \left\{\frac{1}{s_i / b} - \psi'(s_i / b + 1)\right\} \ind_{\{i = j\}} + \frac{1}{b^2} \left[\frac{1}{(1 - \|\bb{s}\|_1) / b} - \psi'\{(1 - \|\bb{s}\|_1) / b + 1\}\right],
\]
where $\psi$ denotes the digamma function and $\psi'$ its derivative.

Note that one has $1/y - \psi'(y + 1) > 0$ for every real $y \in (0, \infty)$ because
\[
\psi'(y + 1) = \int_0^{\infty} \frac{t e^{-y t}}{e^t - 1} \rd t = \int_0^{\infty} \frac{s e^{-s}}{y^2 (e^{s/y} - 1)} \rd s < \int_0^{\infty} \frac{s e^{-s}}{y^2 (s/y)} \rd s = \frac{1}{y},
\]
where the first equality is a consequence of \citet[p.\,260]{MR0167642}, the second equality follows from the change of variable $s = yt$, and the inequality comes from the fact that $e^x - 1 > x$ for every real $x \in (0, \infty)$.

Therefore, the expression for the second-order partial derivatives above shows that the Hessian matrix of $\bb{s}\mapsto \ln \{ K_{{\bb{s}}/{b} + \bb{1}, {(1 - \|\bb{s}\|_1)}/{b} + 1}(\bb{s}) \}$ has the form $a \mathrm{I}_d + c \bb{1}_{d\times d}$, where $a$ and $c$ are strictly positive quantities that depend on $d$, $b$ and $\bb{s}$, $\mathrm{I}_d$ is the identity matrix of order $d$, and $\bb{1}_{d\times d}$ is a $d\times d$ matrix of $1$s. In particular, this means that all the eigenvalues of the Hessian matrix of the map $\bb{s} \mapsto \ln \{ K_{{\bb{s}}/{b} + \bb{1}, {(1 - \|\bb{s}\|_1)}/{b} + 1}(\bb{s}) \}$ are strictly positive, which proves the claim that $\bb{s}\mapsto K_{{\bb{s}}/{b} + \bb{1}, {(1 - \|\bb{s}\|_1)}/{b} + 1}(\bb{s})$ is log-convex.

Now, maximum values of log-convex functions on a convex domain are always attained on the boundary. In the present case, a vector $\bb{s} \in \mathcal{S}_d$ is on the boundary of the simplex when some of its components $s_1, \ldots, s_d$ or $1 - \|\bb{s}\|_1$ are equal to zero. When this happens, the map $\bb{s}\mapsto K_{{\bb{s}}/{b} + \bb{1}, {(1 - \|\bb{s}\|_1)}/{b} + 1}(\bb{s})$ is equal to a lower dimensional version of itself on a lower dimensional simplex where the only relevant components are the ones that are not zero. Then that lower dimensional version is still log-convex because the result in the previous paragraph is valid for every integer $d \in \N$, and hence its maximum values are necessarily attained on the boundary of that lower dimensional simplex.

By iterating the above dimension reduction argument, one deduces that the largest value of the original map $\bb{s}\mapsto K_{{\bb{s}}/{b} + \bb{1}, {(1 - \|\bb{s}\|_1)}/{b} + 1}(\bb{s})$ must necessarily be at one of the $d+1$ corners of the simplex $\mathcal{S}_d$, i.e., the $d+1$ standard basis vectors in $\R^{d+1}$. Indeed, these corners form all the possible boundaries of the $1$-dimensional simplexes in the dimension reduction argument.

Finally, by the symmetry of the map $\bb{s}\mapsto K_{{\bb{s}}/{b} + \bb{1}, {(1 - \|\bb{s}\|_1)}/{b} + 1}(\bb{s})$ in the variables $s_1,\dots,s_d,1 - \|\bb{s}\|_1$, the value of the function is exactly the same at any one of those corners, so that one can choose any one of them to determine the maximum. If  the corner $(\bb{s}^{\star}, 1 - \|\bb{s}^{\star}\|_1) = (0, \ldots, 0, 1)$ is chosen, then one finds
\[
K_{{\bb{s}^{\star}}/{b} + \bb{1}, {(1 - \|\bb{s}^{\star}\|_1)}/{b} + 1}(\bb{s}^{\star}) = \frac{\Gamma(1/b + d + 1)}{\Gamma(1 / b + 1)} = \prod_{i=1}^d (1/b + i).
\]
This completes the proof.

\end{proof}

\vspace{-8mm}
\section*{Acknowledgments}

\noindent
Bertin is supported by grants 1190801 and 1221373 from the Fondo Nacional de Desarrollo Cient\'ifico y Tecnol\'ogico, MATH-AmSud 20MATH05, Proyecto Puente UVA 20993, the Centro de Modelamiento Matem\'atico (ACE210010 and FB210005), and BASAL funds for centers of excellence from Chile's Agencia Nacional de Investigaci\'on y Desarrollo. Genest's research is funded in part by the Canada Research Chairs Program (Grant no.~950--231937), the Natural Sciences and Engineering Research Council of Canada (RGPIN--2016--04720), and the Trottier Institute for Science and Public Policy. Klutchnikoff acknowledges support from the MATH-AmSud 20MATH05 grant. Ouimet benefited from postdoctoral fellowships from the Natural Sciences and Engineering Research Council of Canada and from the Fond qu\'eb\'ecois de la recherche -- Nature et technologies. Ouimet's research is currently funded by a CRM-Simons postdoctoral fellowship from the Centre de recherches math\'ematiques (Montr\'eal) and the Simons Foundation.

%\section*{References}

%\bibliographystyle{authordate1}
\bibliographystyle{myjmva}
\bibliography{BGKO_2022_bib}

%\bibliographystyle{myjmva}
%\section*{}
%or one can use (see the template for details)
%\begin{thebibliography}{99}
%
%\bibitem[{Agresti(2013)}]{Agresti13}
%\bibinfo{author}{A.~Agresti}, \bibinfo{title}{{Categorical Data Analysis}},
%  \bibinfo{publisher}{Wiley, Hoboken}, \bibinfo{year}{2013}.
%%Type = Article
%\bibitem[{Aitchison and Silvey(1958)}]{aitchison1958maximum}
%\bibinfo{author}{J.~Aitchison}, \bibinfo{author}{S.D.~Silvey},
%  \bibinfo{title}{Maximum-likelihood estimation of parameters subject to
%  restraints}, \bibinfo{journal}{The Annals of Mathematical Statistics}
%  \bibinfo{volume}{29} (\bibinfo{year}{1958}) \bibinfo{pages}{813--828}.
% \bibitem{Balsu} A. Balsubramani, S. Dasgupta, Y. Freund, \newblock The fast
%convergence of incremental PCA, \newblock Advances in Neural Information
%Processing Systems 26 (2013) 3174--3182.
%
%\end{thebibliography}
%

\end{document}